\documentclass[12pt]{amsart}
\usepackage{amsmath,amsthm,amssymb}
\usepackage{hyperref, a4wide}
\hypersetup{
    pdftitle=   {Topological classification of torus manifolds which have codimension one extended actions},
    pdfauthor=  {Suyoung Choi and Shintar\^o Kuroki}
}

\theoremstyle{plain}
\newtheorem{theorem}{Theorem}[section]
\newtheorem{proposition}[theorem]{Proposition}
\newtheorem{lemma}[theorem]{Lemma}

\newtheorem{corollary}[theorem]{Corollary}

\newtheorem{problem}[theorem]{Problem}

\theoremstyle{definition}

\newtheorem{example}[theorem]{Example}

\newcommand{\C}{\mathbb{C}}

\newcommand{\Z}{{\mathbb{Z}}}
\newcommand{\Q}{{\mathbb{Q}}}
\newcommand{\R}{{\mathbb{R}}}
\newcommand{\N}{{\mathbb{N}}}

\renewcommand{\iff}{\Leftrightarrow}

\newcommand{\mD}{\mathfrak{D}}
\newcommand{\mE}{\mathfrak{E}}

\newcommand{\mM}{\mathfrak{M}}

\newcommand{\mR}{\mathfrak{R}}

\def\arxiv#1{\href{http://arXiv.org/abs/#1}{arXiv:#1}}

\begin{document}%
\title[Topological classification of some torus manifolds]{Topological classification of torus manifolds which have codimension one extended actions}

\author{Suyoung Choi}
\address{Department of Mathematics, Osaka City University, Sumiyoshi-ku, Osaka 558-8585, Japan}
\email{choisy@sci.osaka-cu.ac.jp}

\author{Shintar\^o Kuroki}
\address{Department of Mathematical Sciences, KAIST, 335 Gwahangno, Yuseong-gu, Daejeon 305-701, Republic of Korea}
\email{kuroki@kaist.ac.kr}

\thanks{The first author was supported by Brain Korea 21 project, KAIST and the Japanese Society for the Promotion of Sciences (JSPS grant no. P09023).
The second author was partially supported by Basic Science Research Program through the National Research Foundation of Korea(NRF) funded by the Ministry of Education, Science and Technology (2010-0001651), Fudan University and the Fujyukai foundation.
}

\keywords{Sphere bundle, Complex projective bundle, Torus manifold, Non-singular toric variety, Quasitoric manifold, Toric topology}
\subjclass[2000]{Primary 57S25; Secondary 55R25}

\date{\today}
\maketitle

\begin{abstract}
A toric manifold is a compact non-singular toric variety equipped with a natural half-dimensional compact torus action.
A torus manifold is an oriented, closed, smooth manifold of dimension $2n$ with an effective action of a compact torus $T^{n}$ having a non-empty fixed point set. Hence, a torus manifold can be thought
of as a generalization of a toric manifold.
In the present paper, we focus on a certain class $\mM$ in the family of torus manifolds with codimension one extended actions, and we give a topological classification of $\mM$. 
As a result, their topological types are completely determined by their cohomology rings and real characteristic classes.

The problem whether the cohomology ring determines the topological type of a toric manifold or not is one of the most interesting open problems in toric topology. 
One can also ask this problem for the class of torus manifolds even if its orbit spaces are highly structured.
Our results provide  a negative answer to this problem for torus manifolds.
However, we find a sub-class of torus manifolds with codimension one extended actions which is not in the class of toric manifolds but which is classified by their cohomology rings.
\end{abstract}

\section{Introduction} \label{sec:1}
A \emph{toric variety} of dimension $n$ is a normal algebraic variety on which an algebraic torus $(\C^\ast)^n$ acts with a dense orbit. A compact non-singular toric variety is called a \emph{toric manifold}. We consider a projective toric manifold and regard the compact torus $T^n$ as the standard compact subgroup in $(\C^\ast)^n$. Then $T^n$ also acts on a projective toric manifold and there is a moment map whose image is a simple convex polytope.
Moreover, one can see that the $T^n$-action is \emph{locally standard},
that is, locally modelled by the standard
$T^n$-action
on $\C^n$.
By taking these two characteristic properties as a starting point, Davis and Januszkiewicz \cite{DJ} first introduced the notion of a quasitoric manifold as topological generalization of a projective toric manifold in algebraic geometry.
A \emph{quasitoric manifold} is a smooth closed manifold of dimension $2n$ with a locally standard
$T^n$-action whose orbit space can be identified with a simple polytope. Note that both toric and quasitoric manifolds have a fixed point.

As an ultimate generalization of both toric and quasitoric manifolds,
Hattori and Masuda \cite{HM} introduced a \emph{torus manifold} (or \emph{unitary toric manifold} in the earlier terminology \cite{M}) which is an oriented, closed, smooth manifold of dimension $2n$ with an effective $T^{n}$-action having a non-empty fixed point set.
The orbit space $M/T^n$ is not necessary a simple polytope.
A compact manifold with corner is \emph{nice} if there are exactly $n$ codimension one faces meeting at each vertex, and is a \emph{homotopy cell} if it is nice and all of its faces are contractible.
A homotopy cell is a natural generalization of a simple polytope.
Therefore, the set of torus manifolds whose orbit space is a homotopy cell is a good family of manifolds for which the toric theory can be developed in the topological category in a nice way, see \cite{MP,MS}.
Obviously both toric or quasitoric manifolds are contained in this family.
%
%It is known by  that the orbit space is a homotopy cell if and only if a torus manifold $M$ satisfies the following two conditions:
%\begin{enumerate}
%  \item $H^{odd}(M) = 0$;
%  \item $M$, $M_i$'s and connected components of any multiple intersection of $M_i$'s are all simply connected,
%\end{enumerate}
%where $M_i$'s are characteristic submanifolds of $M$.

Let $G$ be a compact, connected Lie group with maximal torus $T^{n}$.
In \cite{K,K1,K2}, the second author studied the torus manifolds which have extended $G$-actions.
We call such a torus manifold a {\it torus manifold with $G$-actions}.
There was given
a
complete classification of
the
torus manifolds with $G$-actions whose $G$-orbits have codimension zero or one principal orbits,
up to equivariant diffeomorphism.
If a simply connected torus manifold with $G$-actions has a codimension zero principal orbit, i.e., the $G$-action is transitive,
then such a torus manifold is a product of complex projective spaces and spheres.
This provides also the (non-equivariantly) topological classification.

Let $\widetilde{\mM}$ be the set of simply connected torus manifolds $M^{2n}$
with $G$-actions whose $G$-orbits have codimension one principal orbits.
It follows from results in \cite{K1,K2} that $\widetilde{\mM}$ consists of the following three types of manifolds:
\begin{itemize}
\item TYPE 1: $\prod_{j=1}^{b} S^{2m_j} \times \left( \prod_{i=1}^{a} S^{2\ell_i +1} \times_{(S^{1})^a} P(\C_\rho^{k_1} \oplus \C^{k_2}) \right)$;
\item TYPE 2: $\prod_{j=1}^{b} S^{2m_j} \times \left( \prod_{i=1}^{a} S^{2\ell_i +1} \times_{(S^{1})^a} S(\C_\rho^{k} \oplus \R) \right)$;
\item TYPE 3: $\prod_{j=1}^{b} S^{2m_j} \times \left( \prod_{i=1}^{a} S^{2\ell_i +1} \times_{(S^{1})^a} S(\C_\rho^{k_1} \oplus \R^{2k_2 + 1}) \right)$,
\end{itemize}
where $P(\C_\rho^{k_1} \oplus \C^{k_2})=(\C_\rho^{k_1} \oplus \C^{k_2}-\{0\})/\C^{*}$ is a complex projective space,
$S(\C^{\ell}\oplus \R^{m})\subset \C^{\ell}\oplus \R^{m}$ is a sphere,
and $(S^{1})^{a}$ acts on $\prod_{i=1}^{a} S^{2\ell_i +1}$ naturally
and on $\C_{\rho}^{k_{1}}$ through the group homomorphism $(S^{1})^{a}\to S^{1}$ defined by
\begin{eqnarray*}
(t_{1}, \ldots, t_{a}) \mapsto t_{1}^{\rho_{1}}\cdots t_{a}^{\rho_{a}}
\end{eqnarray*}
for $\rho = (\rho_1, \ldots, \rho_a) \in \Z^a$.

We denote the subset of $\widetilde{\mM}$ satisfying $a=1$ and $b=0$ by $\mM$.
Obviously, a manifold in $\mM$ is a torus manifold and, furthermore, its orbit space is a homotopy cell.
The purpose of this paper is to prove that all manifolds in $\mM$ are classified by their cohomology rings, Pontrjagin classes and Stiefel-Whitney classes up to diffeomorphism (see Theorem~\ref{thm-main}).

We strongly remark that the cohomology ring is not enough to classify them. Recently, the topological classification of toric manifolds has attracted the attention of toric topologists (see \cite{MS}).
Of special interest is the following problem which is now called the \emph{cohomological rigidity problem} for toric manifolds:
\begin{problem}
\label{problem:rigidity}
Are toric manifolds diffeomorphic (or homeomorphic) if their cohomology rings are isomorphic as graded rings?
\end{problem}
One can also ask the problem for quasitoric manifolds and torus manifolds (see \cite[Sections 6-7]{MS}). In the class of toric or quasitoric manifolds, Problem~\ref{problem:rigidity} is still open, apart from some affirmative answers (see e.g. \cite{MP2,CMS2,CM}).

In this paper, we give negative answers to Problem \ref{problem:rigidity} for the class of torus manifolds whose orbit space is a homotopy cell, see Theorem~\ref{thm-main} and Example~\ref{example:counterexample} in Section~\ref{sec:4}.

This paper is organized as follows.
In Section~\ref{sec:2}, we compute cohomology rings and real characteristic classes of
manifolds in $\mM$.
In Section~\ref{sec:3}, we present the main result of this paper.
In Section~\ref{sec:4}, we prove the main result and give the explicit classification of manifolds in $\mM$. We also exhibit several non-trivial examples.
Finally, in Section~\ref{sec:5}, we give a revision of the cohomological rigidity problem for torus manifolds in \cite{MS}.

%%%%%%%%%%%%%%%%%%%%%%%%%%%%%%%%%%%%%%%%%%%%%%%%%%%%%%%%%%%%%%%%%%%%%%%%%%%
%section 2
%%%%%%%%%%%%%%%%%%%%%%%%%%%%%%%%%%%%%%%%%%%%%%%%%%%%%%%%%%%%%%%%%%%%%%%%%%%
\section{Topological invariants}
\label{sec:2}

We will use the following standard symbols in this paper:
\begin{itemize}
\item $H^\ast(X)$ is the cohomology ring of $X$ over $\Z$-coefficients;
\item $w(X)=\sum_{j=0}^\infty w_j(X)$ is the total Stiefel-Whitney class of $X$, where $w_j(X)$ is the $j$-th Stiefel-Whitney class;
\item $p(X)=\sum_{j=0}^\infty p_j(X)$ is the total Pontrjagin class of $X$, where $p_j(X)$ is the $j$-th Pontrjagin class;
\item $\Z[x_{1}, \ldots, x_{m}]$ is the polynomial ring generated by $x_1, \ldots, x_m$;
\item $\langle f_{1}(x_{1}, \ldots, x_{m}),\ \ldots,\ f_{s}(x_{1}, \ldots, x_{m})\rangle$ is the ideal in $\Z[x_{1}, \ldots, x_{m}]$
generated by the polynomials $f_{j}(x_{1}, \ldots, x_{m})$ ($j=1, \ldots, s$);
\item $E(\eta)$ is the total space of a fibre bundle $\eta$.
\end{itemize}

Let $\mM_i$ ($i=1,2,3$) be the subset of $\mM$ of TYPE $i$.
By the definition of $\mM_{i}$, an element $N_i \in \mM_i$ is as follows:
\begin{eqnarray*}
N_{1}&=&S^{2\ell +1} \times_{S^{1}} P(\C_\rho^{k_1} \oplus \C^{k_2}); \\
N_{2}&=&S^{2\ell +1} \times_{S^{1}} S(\C_\rho^{k} \oplus \R); \\
N_{3}&=&S^{2\ell +1} \times_{S^{1}} S(\C_\rho^{k_1} \oplus \R^{2k_2 + 1}),
\end{eqnarray*}
where $\rho \in \Z$.
In this section,
we compute three topological invariants $H^\ast(N_{i})$, $w(N_{i})$ and $p(N_{i})$ of $N_{i}$ for $i=1,\ 2,\ 3$.
%\begin{itemize}
%\item cohomology rings $H^\ast(N_{i})$;
%\item Stiefel-Whiteny classes $w(N_{i})$;
%\item Pontrjagin classes $p(N_{i})$.
%\end{itemize}

%%%%%%%%%%%%%%%%%%%%%%%%%%%%%%%%%%%%%%%%%%%%%%%%%%%%%%%%%%%%%%%%%%%%%%%%%%%
%section 2.1
%%%%%%%%%%%%%%%%%%%%%%%%%%%%%%%%%%%%%%%%%%%%%%%%%%%%%%%%%%%%%%%%%%%%%%%%%%%
\subsection{Topological invariants of $N_{1}$}
\label{sec:2.1}

The purpose of this subsection is to compute topological invariants of
\begin{eqnarray*}
N_{1}=S^{2\ell +1} \times_{S^{1}} P(\C_\rho^{k_1} \oplus \C^{k_2}).
\end{eqnarray*}
In order to compute them, we first recall the torus action on $N_{1}$. Note that, in this case, the dimension of the torus is $\ell+k_{1}+k_{2}-1$.
The torus action on $N_{1}$ is defined as follows ($k_{1}, k_{2}\ge 1$):
\begin{eqnarray*}
& &(a_{1}, \ldots, a_{\ell}, b_{1}, \ldots, b_{k_{1}}, c_{1}, \ldots, c_{k_{2}-1})\cdot[(x_{0}, \ldots, x_{\ell}),
[y_{1}; \ldots; y_{k_{1}}; y_{1}'; \ldots; y_{k_{2}}']] \\
&=&[(x_{0},\ a_{1}x_{1}, \ldots, a_{\ell}x_{\ell}), [b_{1}y_{1}; \cdots;
b_{k_{1}}y_{k_{1}}; c_{1}y_{1}'; \cdots; c_{k_{2}-1}y_{k_{2}-1}'; y_{k_{2}}']],
\end{eqnarray*}
where $a_{i}$, $b_{j}$, $c_{k}\in S^{1}$, $(x_{0}, \ldots, x_{\ell})\in S^{2\ell+1}\subset \C^{\ell+1}$ and
$[y_{1}; \cdots; y_{k_{1}}; y_{1}'; \cdots; y_{k_{2}}']\in P(\C^{k_{1}}_{\rho}\oplus\C^{k_{2}})$.
By this action, we can easily check that $N_{1}$ is a quasitoric manifold
over $\Delta^{\ell}\times \Delta^{k_{1}+k_{2}-1}$ (product of two simplices) whose dimension is $2(\ell+k_{1}+k_{2}-1)$.
Therefore, the {\it Davis-Januszkiewicz formula} \cite[Theorem 4.14, Corollary 6.8]{DJ} can be used to compute topological invariants of quasitoric manifolds.

%%%%%%%%%%%%%%%%%%%%%%%%%%%%%%%%%%%%%%%%%%%%%%%%%%%%%%%%%%%%%%%%%%%%%%%%%%%
%section 2.1.1
%%%%%%%%%%%%%%%%%%%%%%%%%%%%%%%%%%%%%%%%%%%%%%%%%%%%%%%%%%%%%%%%%%%%%%%%%%%
\subsubsection{The Davis-Januszkiewicz formula}
\label{sec:2.1.1}
In this subsection, we review the Davis-Januszkiewicz formula.

Let $P$ be a simple polytope of dimension $n$ with $m$ facets and $M$ a quasitoric manifold of dimension $2n$ over $P$.
As is well-known, the {\it equivariant cohomology ring} $H_{T^{n}}^\ast(M;\ \Z)=H^\ast(ET\times_{T^n}M)$ of $M$ has the following ring structure:
\begin{eqnarray*}
H_{T^{n}}^\ast(M;\ \Z)\simeq \Z[v_{1}, \ldots, v_{m}]/\mathcal{I},
\end{eqnarray*}
where $v_{j}$ ($\deg v_{j}=2$, $j=1, \ldots, m$) is the equivariant Poincar\'e dual
of codimension two invariant submanifold $M_{j}$ in $M^{2n}$ ({\it characteristic submanifolds}) and $\mathcal{I}$ is
an ideal of the polynomial ring $\Z[v_{1}, \ldots, v_{m}]$ generated by $\{\prod_{j\in I}v_{j}\ |\ \bigcap_{j\in I}M_{j}=\emptyset\}$.
This ring $\Z[v_{1}, \ldots, v_{m}]/\mathcal{I}$ is precisely {\it the face ring} of $M/T=P$ by regarding $v_{i}$ as the facet $F_{i}$ in $P$,
because the image of the characteristic submanifold $M_{i}$ of the orbit projection $M\to P$ is the facet $F_{i}$.
Note that $\deg v_{i}=2$.

Let $\pi\colon ET\times_{T}M\to BT$ be the natural projection.
Then, we can define the induced homomorphism
\begin{eqnarray*}
\pi^\ast \colon H^{*}(BT)=\Z[t_{1}, \ldots, t_{n}]\longrightarrow H^\ast(ET\times_{T}M)=H_{T^{n}}^\ast(M^{2n};\ \Z).
\end{eqnarray*}
Moreover, the $\pi^{*}$-image $\pi^\ast(t_{i})$ of $t_{i}$ ($i=1, \ldots, n$) can be described explicitly as follows.
Let $\Lambda$ be an $n\times m$ matrix
$\Lambda=(\lambda_{1}\cdots \lambda_{m})$, where $\lambda_{j}\in \Z^{n}$ ($j=1, \ldots, m$)
corresponds to the generator of Lie algebra of isotropy subgroup of characteristic submanifold $M_{j}$.
We call $\Lambda$ the \emph{characteristic matrix} of $M$.
Put $\lambda_{j}=(\lambda_{1j}, \ldots, \lambda_{nj})^{t}\in \Z^{n}$.
Then we have
\begin{eqnarray}
\label{image of pi^*}
\pi^\ast(t_{i})=\sum_{j=1}^{m}\lambda_{ij}v_{j}.
\end{eqnarray}
Let $\mathcal{J}$ be the ideal in $\Z[v_{1}, \ldots, v_{m}]$ generated by $\pi^{*}(t_{i})$ for all $i=1, \ldots, n$.
Then the ordinary cohomology of quasitoric manifolds has the following ring structure:
\begin{equation}
\label{ordinary-cohom-quasitoric}
H^\ast(M)\simeq \Z[v_{1}, \ldots, v_{m}]/(\mathcal{I}+\mathcal{J}).
\end{equation}
Moreover, for the inclusion $\iota:M\to ET\times_{T}M$,
the Pontrjagin class\footnote{In \cite[Corollary 6.8]{DJ},
the Pontrjagin class of quasitoric manifolds ({\it toric manifolds} in \cite{DJ})
is $\iota^{*}\prod_{i=1}^{m}(1-v_{i}^{2})$.
However, this formula coincides with $1-p_{1}(M)+p_{2}(M)-\cdots=\sum_{i=0}^{m}(-1)^{i}p_{i}(M)$.
Therefore, by \cite{MS},
the Pontrjagin class of quasitoric manifolds must be $p(M)=1+p_{1}(M)+p_{2}(M)+\cdots=\iota^{*}\prod_{i=1}^{m}(1+v_{i}^{2})$}
and the Stiefel-Whitney class
can be described the following {\it Davis-Januszkiewicz formula}:
\begin{equation}
\label{Pontrjagin-class-quasitoric} p(M)=\iota^{*}\prod_{i=1}^{m}(1+v_{i}^{2})\quad \text{ and } \quad w(M)=\iota^{*}\prod_{i=1}^{m}(1+v_{i}).
\end{equation}
%%%%%%%%%%%%%%%%%%%%%%%%%%%%%%%%%%%%%%%%%%%%%%%%%%%%%%%%%%%%%%%%%%%%%%%%%%%
%section 2.1.2
%%%%%%%%%%%%%%%%%%%%%%%%%%%%%%%%%%%%%%%%%%%%%%%%%%%%%%%%%%%%%%%%%%%%%%%%%%%
\subsubsection{Topological invariants of $N_{1}$}
\label{sec:2.1.2}
Now we shall compute $H^\ast(N_1),~p(N_1)$ and $w(N_1)$.
In order to use the Davis-Januszkiewicz formula, we first compute the characteristic matrix of $N_{1}$.

By using $N_{1}/T=\Delta^{\ell}\times \Delta^{k_{1}+k_{2}-1}$ and choosing an appropriate order of its characteristic submanifolds, we may assume that the characteristic matrix of $N_{1}$ is an $n \times (n+2)$ matrix of the form
\begin{equation}\label{ch-1}
\left(
\begin{array}{ccccc}
I_{\ell} & 0 & 0 & \mathbf{1} & 0 \\
0 & I_{k_{1}} & 0 & \rho\mathbf{1} & \mathbf{1} \\
0 & 0 & I_{k_{2}-1} & 0 & \mathbf{1}
\end{array}
\right),
\end{equation}
where $n=\ell+k_{1}+k_{2}-1$ and $\mathbf{1}=(1, \ldots, 1)^{t}\in \Z^r$ for $r=\ell,~k_{1}$ and $k_{2}$.
Since $N_{1}/T=\Delta^{\ell}\times \Delta^{k_{1}+k_{2}-1}$, we can compute its equivariant cohomology as follows:
\begin{equation}
\label{equiv-1}
H_{T^{n}}^{*}(N_{1})\simeq
\Z[v_{1}, \ldots, v_{\ell+1}, w_{1}, \ldots, w_{k_{1}+k_{2}}]/\mathcal{I},
\end{equation}
where $\deg v_{i}=\deg w_{j}=2$ and $\mathcal{I}$ is generated by $v_{1}\cdots v_{\ell+1}$ and $w_{1}\cdots w_{k_{1}+k_{2}}$ (see Section \ref{sec:2.1.1}).
The ideal $\mathcal{J}$ is generated by the following $\ell + k_1+k_2-1$ elements from \eqref{image of pi^*} and \eqref{ch-1}:
\begin{eqnarray}
& &v_{i}+v_{\ell+1}\ (i=1,\ \ldots,\ \ell), \nonumber \\
& &\label{J-1} w_{j}+\rho v_{\ell+1}+w_{k_{1}+k_{2}}\ (j=1,\ \ldots,\ k_{1}), \\
& &w_{j}+w_{k_{1}+k_{2}}\ (j=k_{1}+1,\ \ldots,\ k_{1}+k_{2}-1). \nonumber
\end{eqnarray}
From \eqref{ordinary-cohom-quasitoric}, \eqref{equiv-1} and \eqref{J-1}, we have the ordinary cohomology as follows:
\begin{align*}
H^{*}(N_{1})&\simeq  \Z[v_{1}, \ldots, v_{\ell+1}, w_{1}, \ldots, w_{k_{1}+k_{2}}]/(\mathcal{I}+\mathcal{J}) \\
&=  \Z[v_{\ell+1}, w_{k_{1}+k_{2}}]/ \langle (-1)^{\ell}(v_{\ell+1})^{\ell+1}, (-w_{k_{1}+k_{2}} - \rho v_{\ell+1})^{k_1} (-w_{k_{1}+k_{2}})^{k_2} \rangle \\
&\simeq  \Z[x,y] / \langle x^{\ell +1}, y^{k_2} (y + \rho x)^{k_1} \rangle,
\end{align*}
where $x=v_{\ell+1}$, $y=w_{k_{1}+k_{2}}$.

By \eqref{Pontrjagin-class-quasitoric} and \eqref{J-1},
the real characteristic classes are as follows:
\begin{eqnarray*}
p(N_{1})&=&(1+v_{\ell+1}^{2})^{\ell+1}(1+(\rho v_{\ell+1}+w_{k_{1}+k_{2}})^{2})^{k_{1}}(1+w_{k_{1}+k_{2}}^{2})^{k_{2}} \\
&=& (1+x^{2})^{\ell+1}(1+(\rho x+y)^{2})^{k_{1}}(1+y^{2})^{k_{2}}
\end{eqnarray*}
and
\begin{eqnarray*}
w(N_{1})&\equiv_{2} & (1+v_{\ell+1})^{\ell+1}(1+\rho v_{\ell+1}+w_{k_{1}+k_{2}})^{k_{1}}(1+w_{k_{1}+k_{2}})^{k_{2}} \\
&\equiv_{2} & (1+x)^{\ell+1}(1+\rho x+y)^{k_{1}}(1+y)^{k_{2}}.
\end{eqnarray*}

To summarize, we have the following.
%%%%%%%%%%%%%%%%%%%%%%%%%%%%%%%%%%%%%%%%%%%%%%%%%%%%%%%%%%%%%%%%%%%%%%%%%%%
%Proposition 1
%%%%%%%%%%%%%%%%%%%%%%%%%%%%%%%%%%%%%%%%%%%%%%%%%%%%%%%%%%%%%%%%%%%%%%%%%%%
\begin{proposition}
\label{prop2.1}
Topological invariants of $N_{1}$ are
\begin{align*}
H^{*}(N_1) &= \Z[x,y] / \langle x^{\ell +1}, y^{k_2} (y + \rho x)^{k_1}\rangle, \\
p(N_1) &= (1 + x^2)^{\ell +1}(1+(\rho x + y)^2)^{k_1}(1+y^2)^{k_2} \text{ and} \\
w(N_1) &\equiv_2  (1 + x)^{\ell+1} (1+\rho x + y)^{k_1} (1+y)^{k_2},
\end{align*}
where $\deg x = \deg y = 2$ and $\ell, k_1, k_2 \in \N$.
\end{proposition}

%%%%%%%%%%%%%%%%%%%%%%%%%%%%%%%%%%%%%%%%%%%%%%%%%%%%%%%%%%%%%%%%%%%%%%%%%%%
%section 2.2
%%%%%%%%%%%%%%%%%%%%%%%%%%%%%%%%%%%%%%%%%%%%%%%%%%%%%%%%%%%%%%%%%%%%%%%%%%%
\subsection{Topological invariants of $N_{2}$ and $N_{3}$}
\label{sec:2.2}

The purpose of this subsection is to determine the topological
invariants
of
$$
 N_2=S^{2\ell +1} \times_{S^{1}} S(\C_\rho^{k} \oplus \R)
\quad
\text{ and }
\quad
N_{3}=S^{2\ell +1} \times_{S^{1}} S(\C_\rho^{k_1} \oplus \R^{2k_2 + 1}).
$$
We have the following.
%%%%%%%%%%%%%%%%%%%%%%%%%%%%%%%%%%%%%%%%%%%%%%%%%%%%%%%%%%%%%%%%%%%%%%%%%%%
%Proposition 2
%%%%%%%%%%%%%%%%%%%%%%%%%%%%%%%%%%%%%%%%%%%%%%%%%%%%%%%%%%%%%%%%%%%%%%%%%%%
\begin{proposition}
\label{prop2.2}
Topological invariants of $N_{2}$ are
\begin{eqnarray*}
    H^\ast(N_2)&=& \Z[x,z] / \langle x^{\ell +1}, z(z+ (\rho x)^k )\rangle,\\
    p(N_2) &=& (1 + x^2)^{\ell +1}(1 + \rho^2 x^2)^k \text{ and} \\
    w(N_2) &\equiv_2& (1 + x)^{\ell+1} (1+\rho x )^{k},
\end{eqnarray*}
where $\deg x = 2$, $\deg z = 2k$ and $\ell, k \in \N$.
\end{proposition}

%%%%%%%%%%%%%%%%%%%%%%%%%%%%%%%%%%%%%%%%%%%%%%%%%%%%%%%%%%%%%%%%%%%%%%%%%%%
%Proposition 3
%%%%%%%%%%%%%%%%%%%%%%%%%%%%%%%%%%%%%%%%%%%%%%%%%%%%%%%%%%%%%%%%%%%%%%%%%%%
\begin{proposition}
\label{prop2.3}
Topological invariants of $N_{3}$ are
\begin{eqnarray*}
    H^\ast(N_3)&=& \Z[x,z] / \langle x^{\ell +1}, z^2\rangle, \\
    p(N_3) &=& (1 + x^2)^{\ell +1}(1 +  \rho^2 x^2)^{k_1} \text{ and} \\
    w(N_3) &\equiv_2& (1 + x)^{\ell+1} (1+\rho x )^{k_1},
\end{eqnarray*}
where $\deg x = 2$, $\deg z = 2(k_1 + k_2)$ and $\ell, k_1, k_2 \in \N$.
\end{proposition}

We first prove Proposition \ref{prop2.2}.
Note that $N_{2}$ has the following fibration:
\begin{eqnarray*}
S^{2k}=S(\C_{\rho}^{k}\oplus \R)\longrightarrow N_{2}\stackrel{\pi}{\longrightarrow} S^{2\ell+1}/S^{1}\cong \C P^{\ell}.
\end{eqnarray*}
Thus, $N_{2}$ is a sphere bundle over a complex projective space.
Therefore, we can use the following lemma.
%%%%%%%%%%%%%%%%%%%%%%%%%%%%%%%%%%%%%%%%%%%%%%%%%%%%%%%%%%%%%%%%%%%%%%%%%%%
%Lemma 1
%%%%%%%%%%%%%%%%%%%%%%%%%%%%%%%%%%%%%%%%%%%%%%%%%%%%%%%%%%%%%%%%%%%%%%%%%%%
\begin{lemma}[Duan \cite{Duan}]
\label{lem-Duan}
Let $\pi:E\to M$ be a smooth, oriented $r$-sphere bundle over an oriented manifold $M$
which has a section $s:M\to E$.
Let the normal bundle $\nu$ of the embedding $s$ be oriented by $\pi$, and let
$\chi(\nu)\in H^{r}(M)$ be the Euler class of $\nu$ with respect to this orientation.
Then there exists a unique class $z\in H^{r}(E)$ such that
\begin{eqnarray*}
s^{*}(z)=0\in H^{*}(M)\ and\ <i^{*}(z),\ [S^{r}]>=1.
\end{eqnarray*}
Furthermore,
$H^{*}(E)$ has the basis $\{1,\ z\}$
subject to the relation
\begin{eqnarray*}
z^{2}+\pi^{*}(\chi(\nu))z=0
\end{eqnarray*}
as a module over $H^{*}(M)$.
\end{lemma}

We define a section $s$ of $\pi \colon N_{2}\to \C P^{\ell}$ as
\begin{eqnarray*}
s \colon \C P^{\ell}\ni [z_{0}; \ldots ; z_{\ell}]\mapsto ([z_{0}; \ldots ; z_{\ell}],\ (0, \ldots, 0,\ 1))\in N_{2},
\end{eqnarray*}
where $(0, \ldots, 0,\ 1)\in S(\C_{\rho}^{k}\oplus \R)$.
This map is well-defined because $(0, \ldots, 0,\ 1)\in S(\C_{\rho}^{k}\oplus \R)$ is one of the fixed points of
the $S^{1}$-action on $S(\C_{\rho}^{k}\oplus \R)$.
Then the normal bundle of this section is isomorphic to the following bundle $\xi_{\rho}$:
\begin{eqnarray*}
\C^{k}\longrightarrow S^{2\ell+1}\times_{S^{1}}\C_{\rho}^{k}\longrightarrow \C P^{\ell},
\end{eqnarray*}
where $S^1$ acts on $\C_\rho$ by the representation $t \mapsto t^\rho$.
Then, we have $ \xi_{\rho}\equiv \gamma_{\rho}^{\oplus k}$
where $E(\gamma_{\rho})=S^{2\ell+1}\times_{S^{1}}\C_{\rho}$.
Note that $\gamma_{1}$, i.e., $\rho=1$, is isomorphic to $\overline{\gamma}$ as a complex line bundle, where $\overline{\gamma}$ is the tautological line bundle over $\C P^{\ell}$ with reversed orientation.
Hence, $\gamma_{\rho}=(\overline{\gamma})^{\otimes \rho}$.
%for the orientation reversing, tautological line bundle $\overline{\gamma} (=\gamma_{1})$.
Therefore, the Euler class of $\xi_{\rho}$ is
\begin{eqnarray*}
\chi(\xi_{\rho})=c_{k}(\xi_{\rho})=c_{k}(\gamma_{\rho}^{\oplus k})=c_{1}(\gamma_{\rho})^{k}=c_{1}((\overline{\gamma})^{\otimes \rho})^{k}=(\rho c)^{k}=\rho^{k}c^{k},
\end{eqnarray*}
where $c\in H^{2}(\C P^{\ell})$ is the generator (determined by $c_{1}(\overline{\gamma})$) of the cohomology ring $H^{*}(\C P^{\ell})$.
Using the Gysin exact sequence for the bundle $\pi:N_{2}\to \C P^{\ell}$ and the fact $H^{odd}(\C P^{\ell})=0$
(see e.g. \cite[I-Chapter 3]{MT}),
the induced homomorphism $\pi^\ast$ is injective.
Hence, using Lemma \ref{lem-Duan} and the fact that $H^{*}(\C P^{\ell})=\Z[c]/\langle c^{\ell+1}\rangle$, we have the following two relations in the cohomology ring $H^{*}(N_{2})$:
\begin{eqnarray*}
x^{\ell+1}=0
\quad
\text{ and }
\quad
z^{2}+\rho^{k}x^{k}z=0,
\end{eqnarray*}
for $x=\pi^{*}(c)\in H^{2}(N_{2})$ and some $z\in H^{2k}(N_{2})$.
Making use of the Serre spectral sequence for the bundle $\pi:N_{2}\to \C P^{\ell}$,
there is an epimorphism $\Z[x,\ z]\to H^{*}(N_{2})$, and the cohomology ring of $H^{*}(N_{2})$ coincides with that of
$\C P^{\ell}\times S^{2k}$ as an additive group.
Hence, there is no other relations except those mentioned in the above arguments.
Thus, we have the cohomology ring formula in Proposition~\ref{prop2.2}.

In order to compute characteristic classes,
we regard $N_{2}=S^{2\ell+1} \times_{S^{1}}S(\C_{\rho}^{k}\oplus \R)$
as the unit sphere bundle of the following vector bundle over $\C P^{\ell}$:
\begin{eqnarray}
\label{xi}
\xi=\xi_{\rho}\oplus \underline{\R}\equiv \gamma_{\rho}^{\oplus k}\oplus\underline{\R},
\end{eqnarray}
where $\underline{\R}$ is the trivial
real line bundle. Note that $E(\xi)=S^{2\ell+1}\times_{S^{1}}(\C_{\rho}^{k}\oplus \R)$.
We often denote $N_{2}$ as $S(\xi)$, that is, the unit sphere bundle of $\xi$.

Now $\mathcal{T}$ denotes the tangent bundle of $E(\xi)$.
Then, there is the following pull-back diagram:
\begin{eqnarray*}
\begin{array}{ccc}
\iota^{*}\mathcal{T} & \longrightarrow & \mathcal{T} \\
\downarrow & & \downarrow \\
S(\xi) & \stackrel{\iota}{\longrightarrow} & E(\xi)
\end{array}
\end{eqnarray*}
where $\iota:(N_{2}=)S(\xi)\to E(\xi)$ is the natural inclusion, and we have $\mathcal{T}|_{N_{2}}=\iota^{*}\mathcal{T}=\tau_{2}\oplus \nu_{2}$,
where $\tau_{2}$ is the tangent bundle of $N_{2}=S(\xi)$ and $\nu_{2}$ is the normal bundle of the inclusion $\iota:S(\xi)\to E(\xi)$.
Note that $\nu_{2}$ is a real $1$-dimensional bundle by the equation $\dim E(\xi)-\dim S(\xi)=1$.
Because $N_{2}$ is simply connected, we have the following lemma for $\nu_{2}$ (see \cite{S}).
%%%%%%%%%%%%%%%%%%%%%%%%%%%%%%%%%%%%%%%%%%%%%%%%%%%%%%%%%%%%%%%%%%%%%%%%%%%
%Lemma 2
%%%%%%%%%%%%%%%%%%%%%%%%%%%%%%%%%%%%%%%%%%%%%%%%%%%%%%%%%%%%%%%%%%%%%%%%%%%
\begin{lemma}
\label{lem-trivial}
The vector bundle
$\nu_{2}$ is the trivial real line bundle over $N_{2}$, i.e., $E(\nu_{2})=N_{2}\times \R$.
\end{lemma}

Hence, we have
\begin{eqnarray}
\label{p-tau} \iota^{*}p(\mathcal{T})&=p(\iota^{*}\mathcal{T})=p(\tau_{2}\oplus \nu_{2})=p(\tau_{2})=p(N_{2})
\end{eqnarray}
and
\begin{eqnarray}
\label{w-tau} \iota^{*}w(\mathcal{T})&=w(\iota^{*}\mathcal{T})=w(\tau_{2}\oplus \nu_{2})=w(\tau_{2})=w(N_{2}).
\end{eqnarray}
We have that $\pi:S(\xi)\to \C P^{\ell}$
is decomposed into $\pi=\widetilde{\pi}\circ\iota$ (where $\widetilde{\pi}:E(\xi)\to \C P^{\ell}$) and $\pi^{*}$ is injective; %by the Serre spectral sequence..
therefore, $\iota^{*}:H^{*}(E(\xi))\to H^{*}(S(\xi))$ is injective.
Thus, in order to prove Proposition \ref{prop2.2},
%we compute
it is sufficient to compute
$p(\mathcal{T})$ and $w(\mathcal{T})$.

Let $\widetilde{s}$ be the zero section of $\widetilde{\pi}:E(\xi)\to \C P^{\ell}$.
Consider the following pull-back diagram:
\begin{eqnarray*}
\begin{array}{ccc}
\widetilde{s}^{*}\mathcal{T} & \longrightarrow & \mathcal{T} \\
\downarrow & & \downarrow \\
\C P^{\ell} & \stackrel{\widetilde{s}}{\longrightarrow} & E(\xi).
\end{array}
\end{eqnarray*}
Because the normal bundle $\nu (\C P^{\ell})$ of the image of $\widetilde{s}$ is isomorphic to $\xi$, we have
$$
\widetilde{s}^{*}\mathcal{T}\equiv \tau(\C P^{\ell})\oplus \nu (\C P^{\ell})\equiv \tau(\C P^{\ell})\oplus \xi,
$$ where $\tau(\C P^{\ell})$ is the tangent bundle over $\C P^{\ell}$.
Therefore, by using \eqref{xi} and the well-known real characteristic classes of $\C P^{\ell}$, we have
\begin{align*}
\widetilde{s}^{*}(p(\mathcal{T}))
&=p(\widetilde{s}^{*}\mathcal{T})=p(\tau(\C P^{\ell})\oplus \xi)=p(\tau(\C P^{\ell}))p(\xi) \\
&=(1+c^{2})^{\ell+1}(1+\rho^{2}c^{2})^{k}
\end{align*}
and
\begin{align*}
\widetilde{s}^{*}(w(\mathcal{T}))
&=w(\widetilde{s}^{*}\mathcal{T})=w(\tau(\C P^{\ell})\oplus \xi)=w(\tau(\C P^{\ell}))w(\xi)\\
&=(1+c)^{\ell+1}(1+\rho c)^{k}.
\end{align*}
Because $\widetilde{s}^{*}:H^{*}(E(\xi))\to H^{*}(\C P^{\ell})\simeq \Z[c]/\langle c^{\ell+1}\rangle$ induces the isomorphism
and $\widetilde{s}^{*}=(\widetilde{\pi}^{*})^{-1}$,
we have
\begin{eqnarray*}
p(\mathcal{T})&=(1+\widetilde{\pi}^{*}(c)^{2})^{\ell+1}(1+\rho^{2}\widetilde{\pi}^{*}(c)^{2})^{k}
\end{eqnarray*}
and
\begin{eqnarray*}
w(\mathcal{T})&=(1+\widetilde{\pi}^{*}(c))^{\ell+1}(1+\rho \widetilde{\pi}^{*}(c))^{k}.
\end{eqnarray*}
Hence, since $\iota^{*}\circ\widetilde{\pi}^{*}(c)=\pi^{*}(c)=x$, we have
$$
p(N_{2})=(1+x^{2})^{\ell+1}(1+\rho^{2}x^{2})^{k}
\quad
\text{ and }
\quad
w(N_{2})=(1+x)^{\ell+1}(1+\rho x)^{k}
$$
from \eqref{p-tau} and \eqref{w-tau}.
This establishes Proposition \ref{prop2.2}.

With the method similar to that demonstrated as above, we also have Proposition \ref{prop2.3}.

%%%%%%%%%%%%%%%%%%%%%%%%%%%%%%%%%%%%%%%%%%%%%%%%%%%%%%%%%%%%%%%%%%%%%%%%%%%
%section 3
%%%%%%%%%%%%%%%%%%%%%%%%%%%%%%%%%%%%%%%%%%%%%%%%%%%%%%%%%%%%%%%%%%%%%%%%%%%
\section{Main theorem and Preliminary}
\label{sec:3}

In this section, we state the main theorem and prepare to prove it.

%%%%%%%%%%%%%%%%%%%%%%%%%%%%%%%%%%%%%%%%%%%%%%%%%%%%%%%%%%%%%%%%%%%%%%%%%%%
%section 3.1
%%%%%%%%%%%%%%%%%%%%%%%%%%%%%%%%%%%%%%%%%%%%%%%%%%%%%%%%%%%%%%%%%%%%%%%%%%%
\subsection{Main theorem}
\label{sec:3.1}

In order to state the main theorem, we prepare some notations (also see \cite{P}).
A manifold $M$ in the given family is said to be \emph{cohomologically rigid}
if for any other manifold $M'$ in the family, $H^\ast(M) \simeq H^\ast(M')$ as graded rings implies
that there is
a diffeomorphism $M\cong M'$.
A manifold $M$ in the given family is said to be \emph{rigid by the cohomology ring and the Pontrjagin class}
(resp. \emph{the Stiefel-Whitney class})
if for any other manifold $M'$ in the family, the ring isomorphism $\phi:H^{*}(M; \Z) \simeq H^{*}(M'; \Z)$
such that $\phi(p(M))=p(M')$ (resp.\ $\phi(w(M))=w(M')$) implies
that there is
a diffeomorphism $M\cong M'$.
We remark that if $M$ is cohomologically rigid in the given family, then $M$ is automatically
rigid by the cohomology ring and the Pontrjagin class (and the Stiefel-Whitney class).
Throughout this
and the next
sections,
we put
\begin{eqnarray*}
A(\ell,\rho,k_{1},k_{2})=S^{2\ell+1}\times_{S^{1}}P(\C^{k_{1}}_{\rho}\oplus\C^{k_{2}})\in\mM_{1}
\end{eqnarray*}
for some $\rho\in \Z$, $k_{1}\in \N$ and $k_{2}\in \N$,
and
\begin{eqnarray*}
B(\ell,\rho,k_{1},k_{2})=S^{2\ell+1}\times_{S^{1}}S(\C_{\rho}^{k_{1}}\oplus\R^{2k_{2}+1})\in\mM_{2}\cup\mM_{3},
\end{eqnarray*}
for some $\rho\in \Z$, $k_{1}\in \N$ and non-negative integer $k_{2}(\ge 0$).

Moreover, let the subsets $\mR_{1}$, $\mR_{2}$ and $\mR_{3}$ in $\mM$ be defined as follows:
\begin{itemize}
\item $\mR_{1}$ is the set of manifolds which are cohomologically rigid;
\item $\mR_{2}$ is the set of manifolds which are not cohomologically rigid, but rigid by the cohomology ring and the Pontrjagin class;
\item $\mR_{3}$ is the set of manifolds which are not rigid by the cohomology ring and the Pontrjagin class, but rigid by the cohomology ring and the Stiefel-Whitney class in $\mM$.
\end{itemize}
Note that $\mR_{i}\cap \mR_{j}=\emptyset$ for $i\not=j$.

Now we may state the main theorem.
%%%%%%%%%%%%%%%%%%%%%%%%%%%%%%%%%%%%%%%%%%%%%%%%%%%%%%%%%%%%%%%%%%%%%%%%%%%
%Main theorem
%%%%%%%%%%%%%%%%%%%%%%%%%%%%%%%%%%%%%%%%%%%%%%%%%%%%%%%%%%%%%%%%%%%%%%%%%%%

\begin{theorem}
\label{thm-main}
Diffeomorphism types of $\mM$ are completely determined by their cohomology rings, Pontrjagin classes and Stiefel-Whitney classes.

In particular, $\mM=\mR_{1}\sqcup \mR_{2}\sqcup \mR_{3}$ and
\begin{enumerate}
\item the subset $\mR_{1}$ consists of the following manifolds:
\begin{eqnarray*}
& & A(\ell,\rho,k_{1},k_{2}); \\
& & B(\ell,\rho,1,0); \\
& & B(\ell,\rho,k_{1},0) \quad {\rm for}\ \rho\not=0\ {\rm and}\ 4\le 2k_{1}\le \ell,
\end{eqnarray*}
\item the subset $\mR_{2}$ consists of the following manifolds:
\begin{eqnarray*}
& & B(\ell,\rho,k_{1},0)\quad {\rm for}\ \rho\not=0\ {\rm and}\ 3\le\ell+1\le 2k_{1}; \\
& & B(\ell,0,k_{1},0)\quad {\rm for}\ \ell\ge 2\ {\rm and}\ k_{1}\ge 2; \\
& & B(\ell,\rho,k_{1},k_{2})\quad {\rm for}\ \ell\ge 2\ {\rm and}\ k_{2}>0,
\end{eqnarray*}
\item the subset $\mR_{3}$ consists of the following manifolds:
\begin{eqnarray*}
B(1,\rho,k_{1},k_{2})\quad {\rm for}\  k_{1}+k_{2}\ge 2.
\end{eqnarray*}
\end{enumerate}
\end{theorem}

We also have an explicit topological classification of $\mM$. See Corollaries \ref{CASE1(top)}, \ref{top classification of CASE2(1)}, \ref{top classification of CASE2(2)} and \ref{top classification of CASE2(3)} in Section~\ref{sec:4}.
Due to Theorem \ref{thm-main}, we have torus manifolds which do not satisfy the cohomological rigidity even if their orbit spaces are homotopy cells. They provide
negative answers
to the cohomological rigidity problem of torus manifolds,
(a problem which first appeared in \cite[Problem 1 and Section 7]{MS}).

%%%%%%%%%%%%%%%%%%%%%%%%%%%%%%%%%%%%%%%%%%%%%%%%%%%%%%%%%%%%%%%%%%%%%%%%%%%
%section 3.2
%%%%%%%%%%%%%%%%%%%%%%%%%%%%%%%%%%%%%%%%%%%%%%%%%%%%%%%%%%%%%%%%%%%%%%%%%%%
\subsection{Preliminaries}
\label{sec:3.2}

In this subsection, we prepare to prove Theorem \ref{thm-main}.

Due to the definition of $N_1\in \mM_{1}$,
the manifold $N_{1}$ is the projectivization of the following vector bundle $\eta_{\rho}$:
\begin{eqnarray*}
\C^{k_{1}+k_{2}}\longrightarrow S^{2\ell +1} \times_{S^{1}} (\C_\rho^{k_1} \oplus \C^{k_2})\longrightarrow \C P^{\ell}.
\end{eqnarray*}
Now we have
\begin{eqnarray*}
\eta_{\rho}\equiv (\gamma^{\otimes (-\rho)})^{\oplus k_1} \oplus \underline{\C}^{\oplus k_2}
\end{eqnarray*}
 over $\C P^{\ell}$ where $\gamma$ is the tautological line bundle, i.e.,
$E(\gamma)=S^{2\ell+1}\times_{S^{1}}\C_{(-1)}$,
and $\underline{\C}$ is the trivial complex line bundle.
Thus, $\mM_1$ consists of {\it 2-stage generalized Bott manifolds}, i.e., projectivizations of
Whitney sums of line bundles over complex projective spaces (see \cite{CMS,CMS2}).
Therefore, we can use the following theorem.
%%%%%%%%%%%%%%%%%%%%%%%%%%%%%%%%%%%%%%%%%%%%%%%%%%%%%%%%%%%%%%%%%%%%%%%%%%%
%Theorem 3.2
%%%%%%%%%%%%%%%%%%%%%%%%%%%%%%%%%%%%%%%%%%%%%%%%%%%%%%%%%%%%%%%%%%%%%%%%%%%
\begin{theorem}[Choi-Masuda-Suh \cite{CMS2}]
\label{thm for 1}
2-stage generalized Bott manifolds are diffeomorphic if and only if their integral cohomology rings are isomorphic.
\end{theorem}

In order to show whether the manifold in $\mM$ is cohomologically rigid or not, we first compare the cohomology ring of each manifold in $\mM$.
Now, for each $M \in \mM$, its cohomology ring is
\begin{eqnarray*}
H^\ast(M) = \Z[x, w] / \langle x^{\ell+1}, f(x, w)\rangle,
\end{eqnarray*}
where $f$ is a homogeneous polynomial and $\deg x = 2$ and $w=y$ in the case of Proposition \ref{prop2.1} or $w=z$ in the case of
Propositions \ref{prop2.2}, \ref{prop2.3}, that is,
$f$ is one of the following:
\begin{description}
\item[Proposition \ref{prop2.1}] $f(x,\ y)=y^{k_{2}}(y+\rho x)^{k_{1}}$, $\deg y=2$ for $k_{1},\ k_{2}\in \N$;
\item[Proposition \ref{prop2.2}] $f(x,\ z)=z(z+(\rho x)^{k})$, $\deg z=2k$ for $k\in \N$;
\item[Proposition \ref{prop2.3}] $f(x,\ z)=z^{2}$, $\deg z=2k_{1}+2k_{2}$ for $k_{1},\ k_{2}\in \N$.
\end{description}
It easily follows from this fact that
if $H^{*}(M)\simeq H^{*}(M')$ for $M'\in \mM$ such that $H^{*}(M')=\Z[x', w'] / \langle (x')^{\ell'+1}, f'(x', w')\rangle$
then $\deg w=\deg w' (\ge 2)$ and $\deg f=\deg f' (\ge 4)$;
moreover, if $\deg w=\deg w'> 2$ then $\ell=\ell'\ge 1$.
Thus, we may divide the proof of Theorem \ref{thm-main} into the following two cases by the degree of $w$:
\begin{description}
\item[CASE 1] $\deg w=2$, i.e., $2$-dimensional sphere bundle or complex projective bundle;
\item[CASE 2] $\deg w>2$, i.e., $m$-dimensional sphere bundle and $m=\deg w>2$.
\end{description}
Moreover, we divide CASE 2 into the following three sub-cases by $\ell$:
\begin{description}
\item[CASE 2 (1)] $\deg w >2$ and $\ell \ge 4$;
\item[CASE 2 (2)] $\deg w >2$ and $\ell =2,\ 3$;
\item[CASE 2 (3)] $\deg w >2$ and $\ell=1$.
\end{description}

%%%%%%%%%%%%%%%%%%%%%%%%%%%%%%%%%%%%%%%%%%%%%%%%%%%%%%%%%%%%%%%%%%%%%%%%%%%
%section 4
%%%%%%%%%%%%%%%%%%%%%%%%%%%%%%%%%%%%%%%%%%%%%%%%%%%%%%%%%%%%%%%%%%%%%%%%%%%
\section{Proof of the main theorem}
\label{sec:4}

In this section, we prove Theorem \ref{thm-main}.

%%%%%%%%%%%%%%%%%%%%%%%%%%%%%%%%%%%%%%%%%%%%%%%%%%%%%%%%%%%%%%%%%%%%%%%%%%%
%section 4.1
%%%%%%%%%%%%%%%%%%%%%%%%%%%%%%%%%%%%%%%%%%%%%%%%%%%%%%%%%%%%%%%%%%%%%%%%%%%
\subsection{CASE 1 : $\deg w=2$}
\label{sec:4.1}

In this subsection, we assume $\deg w=2$. We first prove that CASE 1 satisfies the cohomological rigidity.
Because $\deg w=2$, this case is a $2$-dimensional sphere bundle or a complex projective bundle over $\C P^{\ell}$ , i.e.,
\begin{eqnarray*}
N_{1}=A(\ell,\rho,k_{1},k_{2})=S^{2\ell+1}\times_{S^{1}}P(\C^{k_{1}}_{\rho}\oplus\C^{k_{2}})
\end{eqnarray*}
or
\begin{eqnarray*}
N_{2}'=B(\ell,\rho,1,0)=S^{2\ell+1}\times_{S^{1}}S(\C_{\rho}\oplus\R).
\end{eqnarray*}
Analyzing the torus action on $N_{2}'$, one can easily prove that $N_{2}'$ is a quasitoric manifold.
Moreover, we see that the orbit space and the characteristic function of $N_{2}'$ coincide with those of $N_{1}$ with $k_{1}=1=k_{2}$, i.e.,
$A(\ell,\rho,1,1)$.
Recall that $N_{1}\in \mM_{1}$ is a $2$-stage generalized Bott manifold (see Section \ref{sec:3.2}).
Hence, by using the construction of Davis-Januszkiewicz in \cite{DJ}, a quasitoric manifold $N_{2}'$ is equivariantly homeomorphic to a $2$-stage generalized Bott manifold $A(\ell,\rho,1,1)$.

We claim that $N_{1}=A(\ell,\rho,1,1)$ and $N_{2}'=B(\ell,\rho,1,0)$ are diffeomorphic.
Let $P_{(1,-1)}(\C_{\rho}\oplus\C)$ be the weighted projective space induced by the identity
$(z_{1},z_{2})\sim (\lambda z_{1},\lambda^{-1}z_{2})$, where $(z_{1},z_{2})\in S(\C_{\rho}\oplus\C)=S^{3}\subset \C_{\rho}\oplus\C$ and $\lambda\in S^{1}$.
Then one can easily see that $P_{(1,-1)}(\C_{\rho}\oplus\C)$ and $S(\C_{\rho}\oplus\R)$ are equivariantly diffeomorphic; therefore, we identify $P_{(1,-1)}(\C_{\rho}\oplus\C)$ and $S(\C_{\rho}\oplus\R)$.
Let us consider the following diffeomorphism
\begin{eqnarray*}
S^{2\ell+1}\times_{S^{1}} S(\C_{\rho}\oplus\C) \stackrel{\psi}{\longrightarrow} S^{2\ell+1}\times_{S^{1}} S(\C_{\rho}\oplus\C)
\end{eqnarray*}
such that
\begin{eqnarray*}
\psi[X,(z_{1},z_{2})]=[X,(z_{1},\overline{z_{2}})],
\end{eqnarray*}
where $X\in S^{2\ell+1}$ and $(z_{1},z_{2})\in S(\C_{\rho}\oplus\C)=S^{3}\subset \C_{\rho}\oplus\C$.
Then we have the induced diffeomorphism $\overline{\psi}:N_{1}\to N_{2}'$ by the following commutative diagram:
\begin{eqnarray*}
\begin{array}{ccc}
S^{2\ell+1}\times_{S^{1}} S(\C_{\rho}\oplus\C) & \stackrel{\psi}{\longrightarrow} & S^{2\ell+1}\times_{S^{1}} S(\C_{\rho}\oplus\C) \\
\pi \downarrow & & \downarrow \pi' \\
A(\ell,\rho,1,1) & \stackrel{\overline{\psi}}{\longrightarrow} & S^{2\ell+1}\times_{S^{1}} P_{(1,-1)}(\C_{\rho}\oplus\C)
\end{array}
\end{eqnarray*}
where $\pi$ is the quotient of fibre $S(\C_{\rho}\oplus\C)$ by the standard multiplication of $S^{1}$, and
$\pi'$ is that by the multiplication of $S^{1}$ induced from the definition of $P_{(1,-1)}(\C_{\rho}\oplus\C)$.
Therefore, we have that $N_{2}'$ is diffeomorphic to a $2$-stage generalized Bott manifold,
in other words, $N_{2}'\in \mM_{1}\cap \mM_{2}$.

Hence, by using Theorem \ref{thm for 1}, the following proposition holds.
%%%%%%%%%%%%%%%%%%%%%%%%%%%%%%%%%%%%%%%%%%%%%%%%%%%%%%%%%%%%%%%%%%%%%%%%%%%
%Proposition 4
%%%%%%%%%%%%%%%%%%%%%%%%%%%%%%%%%%%%%%%%%%%%%%%%%%%%%%%%%%%%%%%%%%%%%%%%%%%
\begin{proposition}
\label{CASE1}
If the cohomology ring of $M$ satisfies
\begin{eqnarray*}
H^\ast(M) = \Z[x, w] / \langle x^{\ell+1}, w^{k_{2}}(w+\rho x)^{k_{1}}\rangle
\end{eqnarray*}
such that $\deg x=\deg w=2$,
then $M$ is diffeomorphic to one of the elements in $\mM_{1}$.
Furthermore, such manifold $M$ is cohomologically rigid in $\mM$, i.e., $M\in\mR_{1}$.
\end{proposition}

We next give the explicit topological classification of CASE 1.
By Proposition \ref{CASE1}, if $M\in \mM$ satisfies CASE 1 then we may put
\begin{eqnarray*}
M\cong S^{2\ell+1}\times_{S^{1}}P(\C^{k_{1}}_{\rho}\oplus\C^{k_{2}})= A(\ell,\rho,k_1,k_2).
\end{eqnarray*}
Therefore, by \cite[Theorem 6.1]{CMS2}, we have the following corollary.

%%%%%%%%%%%%%%%%%%%%%%%%%%%%%%%%%%%%%%%%%%%%%%%%%%%%%%%%%%%%%%%%%%%%%%%%%%%
%Corollary 1
%%%%%%%%%%%%%%%%%%%%%%%%%%%%%%%%%%%%%%%%%%%%%%%%%%%%%%%%%%%%%%%%%%%%%%%%%%%
\begin{corollary}
\label{CASE1(top)}
For the element in $\mM_{1}$,
the following two statements hold:
\begin{enumerate}
\item
The manifold $A(\ell,\rho,k_1,k_2)\cong A(\ell,\rho,k_1,k_2) $ is decomposable, i.e.,
$A(\ell,\rho,k_1,k_2) $ is diffeomorphic to $\C P^\ell \times \C P^{k_1 + k_2}$ if and only if
\begin{eqnarray*}
\left\{
  \begin{array}{ll}
    \rho = 0, & {\rm if}\ \ell>1; \\
    \rho \equiv 0 \mod k_1+k_2+1, & {\rm if}\ \ell=1.
  \end{array}
\right.
\end{eqnarray*}
\item
If $A(\ell,\rho,k_1,k_2)$ is indecomposable, then
$A(\ell,\rho,k_1,k_2)$ is diffeomorphic to $A(\ell',\rho',k_1',k_2')$ if and only if
$\ell=\ell'$, $k_1 + k_2 = k_1' + k_2'$ and there exist $\epsilon=\pm1$ and $r \in \Z$ such that
\begin{eqnarray*}
    (1+ \epsilon r x)^{1+k_2} (1+\epsilon (\rho + r)x)^{k_1} = (1 + \rho'x)^{k_1'} \in \Z[x]/x^{\ell+1}.
\end{eqnarray*}
\end{enumerate}
\end{corollary}

In particular, if $\ell=1$, then the equation in the second statement of Corollary \ref{CASE1(top)} is
equivalent to $\rho k_1 \equiv \pm  \rho' k_1' \text{ (mod $k_1 + k_2 + 1$)}$.
%%%%%%%%%%%%%%%%%%%%%%%%%%%%%%%%%%%%%%%%%%%%%%%%%%%%%%%%%%%%%%%%%%%%%%%%%%%
%section 4.2
%%%%%%%%%%%%%%%%%%%%%%%%%%%%%%%%%%%%%%%%%%%%%%%%%%%%%%%%%%%%%%%%%%%%%%%%%%%
\subsection{CASE 2: $\deg w >2$}
\label{sec:4.2}

In this subsection, we prepare to analyze CASE 2 (1)--(3).
Henceforth, we assume $\deg w >2$.

In CASE 2, a torus manifold $M$ is a manifold in $\mM_2\backslash (\mM_{1}\cap \mM_{2})$ or $\mM_3$, i.e., $M$ is in $\mM_{1}^{c}$ (the complement of $\mM_{1}$) because $\mM_{1}\cap \mM_{3}=\emptyset$.
Moreover, without loss of generality, we may put $M$ as follows:
\begin{eqnarray*}
M= S((\gamma^{\otimes (-\rho)})^{\oplus k_{1}} \oplus \R^{2k_{2}+1})=B(\ell,\rho,k_{1},k_{2}).
\end{eqnarray*}
Because $\deg w=2k_{1}+2k_{2}>2$, we have $(k_{1},\ k_{2})\not=(1,\ 0)$.

Let $M_1=B(\ell_{1},\rho_{1},k_{11},k_{12})$ and $M_2=B(\ell_{2},\rho_{2},k_{21},k_{22})$ be two manifolds in $\mM_{1}^{c}$.
In order to analyze the rigidities in $\mM_{1}^{c}$, we assume $H^{*}(M_{1})\simeq H^{*}(M_{2})$.
Now the cohomology ring $H^{*}(M_{i})$ ($i=1,\ 2$) is given by
\begin{eqnarray}
\label{cohomology}
H^{*}(M_{i})= \Z[x_i, w_i] / \langle x_i^{\ell_{i}+1}, f_i(x_i, w_i)\rangle
\end{eqnarray}
where $\deg x_{i}=2$, $\deg w_{i}=2k_{i1}+2k_{i2}>2$ (see Section \ref{sec:3.2}).
Because $M_{1}$ and $M_{2}$ are sphere bundles over complex projective spaces and $H^{*}(M_{1})\simeq H^{*}(M_{2})$,
one can easily prove that $\ell_{1}=\ell_{2}$ (the dimension of the base space) and
$\deg w_{1}=2k_{11}+2k_{12}=2k_{21}+2k_{22}=\deg w_{2}$ (the dimension of the fibre).
Put $\ell=\ell_{1}=\ell_{2}$.
Then we have
\begin{eqnarray*}
& &H^{*}(M_{1})= \Z[x_1, w_1] / \langle x_1^{\ell+1}, f_1(x_1, w_1)\rangle \\
&\simeq & H^{*}(M_{2})=\Z[x_2, w_2] / \langle x_2^{\ell+1}, f_2(x_2, w_2)\rangle.
\end{eqnarray*}
Here, we may assume that the polynomial $f_{i}$ ($i=1,\ 2$) is
\begin{eqnarray}
\label{f_i}
f_{i}(x_{i},\ w_{i}) = \left\{
                             \begin{array}{lr}
                               w_{i}(w_{i}+(\rho_{i}x_{i})^{k_{i1}}), & \hbox{ if $k_{i2}=0$;} \\
                               w_{i}^2, & \hbox{if $k_{i2}>0$.}
                             \end{array}
                           \right.
\end{eqnarray}

Let $\phi : H^{*}(M_1) \rightarrow H^{*}(M_2)$ be a ring isomorphism.
Because
$H^{2}(M_{1})\stackrel{\phi}{\simeq }H^{2}(M_{2})$ and $\deg w_{i}>2$,
we have
\begin{eqnarray}
\label{eq for x}
\phi(x_1) = \pm x_2.
\end{eqnarray}
If $\phi$ preserves the Pontrjagin classes, then we have
\begin{align}
\label{p(M_i)}
\phi(p(M_{1})) & = (1+x_{2}^2)^{\ell+1}(1+\rho_{1}^{2}x_{2}^{2})^{k_{11}} \nonumber \\
    &=(1+x_{2}^2)^{\ell+1}(1+\rho_{2}^{2}x_{2}^{2})^{k_{21}} = p(M_2) \text{ in $H^\ast(M_2)$}
\end{align}
by Propositions \ref{prop2.2}, \ref{prop2.3} and \eqref{eq for x}.

%%%%%%%%%%%%%%%%%%%%%%%%%%%%%%%%%%%%%%%%%%%%%%%%%%%%%%%%%%%%%%%%%%%%%%%%%%%
%section 4.3
%%%%%%%%%%%%%%%%%%%%%%%%%%%%%%%%%%%%%%%%%%%%%%%%%%%%%%%%%%%%%%%%%%%%%%%%%%%
\subsection{CASE 2 (1) : $\deg w >2$ and $\ell \ge 4$}
\label{sec:4.3}

Assume $\ell\ge 4$.
In this case, we have $x_2^4 \neq 0$ in $H^\ast(M_2)$ because $\ell \geq 4$.

We first prove that this case is always rigid by the cohomology ring and the Pontrjagin class.
Assume that the cohomology ring isomorphism $\phi$ preserves the Pontrjagin classes of $M_{1}$ and $M_{2}$.
Using \eqref{p(M_i)} and $x_2^4 \neq 0$, we have
$$
\phi( p_1(M_1)) = p_1(M_2) \iff k_{11} \rho_1^2 =k_{21} \rho_2^2
$$
and
$$
\phi( p_2(M_1)) = p_2(M_2) \iff {k_{11} \choose 2 } \rho_1^4 = {k_{21} \choose 2 } \rho_2^4.
$$
Therefore, since $k_{11}, k_{21} \in \N$, %and $\phi$ preserves the Pontrjagin classes,
one can easily show that we have either
\begin{eqnarray*}
\rho_1 = \rho_2 = 0
\end{eqnarray*}
 or
\begin{eqnarray*}
\rho_1 = \pm \rho_2\not=0,\quad k_{11}=k_{21}.
\end{eqnarray*}
In both cases, the vector bundle $\left(\gamma^{\otimes (-\rho_1)}\right)^{\oplus k_{11}} \oplus \R^{2k_{12}+1}$
and $\left(\gamma^{\otimes (-\rho_2)}\right)^{\oplus k_{21}} \oplus \R^{2k_{22}+1}$ are isomorphic as a real vector bundle.
This implies that $M_1$ and $M_2$ (unit sphere bundles of these vector bundles) are diffeomorphic.
Note that if $f\colon M_1 \to M_2$ is a diffeomorphism, then its induced isomorphism $f^\ast \colon H^\ast(M_2) \to H^\ast(M_1)$ must preserve Pontrjagin classes.
Hence, we have the following lemma.
%%%%%%%%%%%%%%%%%%%%%%%%%%%%%%%%%%%%%%%%%%%%%%%%%%%%%%%%%%%%%%%%%%%%%%%%%%%
%Lemma 3
%%%%%%%%%%%%%%%%%%%%%%%%%%%%%%%%%%%%%%%%%%%%%%%%%%%%%%%%%%%%%%%%%%%%%%%%%%%
\begin{lemma}
\label{rigid by C-P}
Assume $\deg w>2$ and $\ell\ge 4$.
Then,
there is a graded ring isomorphism $\phi:H^{*}(M_{1})\to H^{*}(M_{2})$ such that $\phi(p(M_{1}))=p(M_{2})$
if and only if
$M_1$ and $M_2$ are diffeomorphic.
\end{lemma}

The above arguments together with Lemma~\ref{rigid by C-P} also provide the explicit topological classification of CASE 2 (1).
%%%%%%%%%%%%%%%%%%%%%%%%%%%%%%%%%%%%%%%%%%%%%%%%%%%%%%%%%%%%%%%%%%%%%%%%%%%
%Corollary 2
%%%%%%%%%%%%%%%%%%%%%%%%%%%%%%%%%%%%%%%%%%%%%%%%%%%%%%%%%%%%%%%%%%%%%%%%%%%
\begin{corollary}
\label{top classification of CASE2(1)}
For elements in $\mM_{1}^{c}$ with $\deg w>2$ and $\ell\ge 4$, we have that
$B(\ell,\rho,k_{1},k_{2})\cong B(\ell',\rho',k_{1}',k_{2}')$ if and only if
$\ell=\ell'$ and one of the following is satisfied.
\begin{enumerate}
\item $\rho=\rho'=0$ and $k_{1}+k_{2}=k_{1}'+k_{2}'$;
\item $\rho=\pm\rho'\not=0$, $k_{1}=k_{1}'$ and $k_{2}=k_{2}'$.
\end{enumerate}
\end{corollary}

In order to check whether the cohomological rigidity holds or not,
we divide this case into the following two cases.

%%%%%%%%%%%%%%%%%%%%%%%%%%%%%%%%%%%%%%%%%%%%%%%%%%%%%%%%%%%%%%%%%%%%%%%%%%%
%section 4.3.1
%%%%%%%%%%%%%%%%%%%%%%%%%%%%%%%%%%%%%%%%%%%%%%%%%%%%%%%%%%%%%%%%%%%%%%%%%%%
\subsubsection{The case when $M$ satisfies $H^{*}(M)\simeq H^{*}(\C P^{\ell}\times S^{2\deg w})$}
\label{sec:4.3.1}
First, we find all manifolds $M\in\mM$ satisfying $H^{*}(M)\simeq H^{*}(\C P^{\ell}\times S^{2\deg w})$.
The symbol $\mD$ represents the subset of such manifolds in $\mM_{1}^{c}$.

Using \eqref{cohomology} and \eqref{f_i}, it is easy to show that the following three cases are in $\mD$:
\begin{enumerate}
\item $M=B(\ell,0,k_{1},k_{2})\cong \C P^{\ell}\times S^{2k_{1}+2k_{2}}$;
\item $M=B(\ell,\rho,k_{1},k_{2})$ with $\rho\not=0$ and $k_{2}>0$;
\item $M=B(\ell,\rho,k_{1},0)$ with $\rho\not=0$ and $k_{1}\ge \ell +1$.
\end{enumerate}

Assume that $M\in \mD$ does not belong in the above three cases.
Then, we may assume that
$$
M=B(\ell,\rho,k_{1},0)\quad \text{ for $k_{1}\le \ell, \rho\not=0$}
$$
such that
\begin{eqnarray*}
H^{*}(M_{1})&=&\Z[x_{1},w_{1}]/\langle x_{1}^{\ell+1},w_{1}(w_{1}+\rho^{k_{1}}x_{1}^{k_{1}})\rangle \\
&\stackrel{\phi}{\simeq} & \Z[x_{2},w_{2}]/\langle x_{2}^{\ell+1},w_{2}^{2}\rangle,
\end{eqnarray*}
where $\deg x_{i}=2$ and $\deg w_{i}=2k_{1}$ for $i=1,\ 2$.
By \eqref{eq for x}, we may put
$$
\phi(x_{1})=\pm x_{2} \quad \text{ and } \quad \phi(w_{1})=a x_{2}^{k_{1}}+b w_{2},
$$
for some $a,\ b\in \Z$.
Because $\phi(w_{1})$ and $x_{2}^{k_{1}}$
must be generators in $H^{2k_{1}}(M_{2})$,
one can see that $b= \pm1 \neq 0$. Because $w_{1}(w_{1}+\rho^{k_{1}}x_{1}^{k_{1}})\equiv 0$, $x_{2}^{\ell+1}\equiv 0$ and $w_{2}^{2}\equiv 0$,
we have that
\begin{align*}
\phi(w_{1}(w_{1}+\rho^{k_{1}}x_{1}^{k_{1}})) &= (a x_{2}^{k_{1}}+b w_{2})(a x_{2}^{k_{1}} +b w_{2}+\rho^{k_{1}}(\pm x_{2})^{k_{1}})  \\
&\equiv a (a+(\pm \rho)^{k_{1}})x_{2}^{2k_{1}} +b (2a+(\pm\rho)^{k_{1}}) w_{2}x_{2}^{k_{1}}\\
&\equiv 0.
\end{align*}
Since $w_{2}x_{2}^{k_{1}}\not=0$, we have $a= - \frac{(\pm\rho)^{k_{1}}}{2}\in \Z$.
Note that $a(a+ (\pm \rho)^{k_1})\not=0$ because $\rho \neq 0$. Hence, we also have $x_{2}^{2k_{1}}=0$.
It follows from the ring structure of $H^{*}(M_{1})$ that $2k_{1}\ge \ell+1$.
Therefore, if $\rho\equiv_{2}0$, $\rho\not=0$ and $k_{1}<\ell+1\le 2k_{1}$,
then the cohomology ring of $M=B(\ell,\rho,k_{1},0)$ and that of $\C P^{\ell}\times S^{2k_{1}}$ are the same.

In summary, if $\deg w>2$ and $\ell\ge 4$, then the cohomology ring $H^{*}(M)$ is isomorphic to $H^{*}(\C P^{\ell}\times S^{2k_{1}+2k_{2}})$ if and only if
$M$ satisfies
one of the following:
\begin{enumerate}
\item $M=B(\ell,0,k_{1},k_{2})\cong \C P^{\ell}\times S^{2k_{1}+2k_{2}}$;
\item $M=B(\ell,\rho,k_{1},k_{2})$ with $\rho\not=0$ and $k_{2}>0$;
\item $M=B(\ell,\rho,k_{1},0)$ with $\rho\not=0$ and $k_{1}\ge \ell +1$;
\item $M=B(\ell,\rho,k_{1},0)$ with $\rho\not=0$, $\rho\equiv_{2}0$ and $k_{1}<\ell+1\le 2k_{1}$.
\end{enumerate}
In other words, $\mD\subset \mM_{1}^{c}$ coincides with the set of manifolds which satisfy 1--4 as above.
Note that the cohomology ring of $M=B(\ell,\rho,k_{1},k_{2})\in\mD$ is isomorphic to $H^{*}(\C P^{\ell}\times S^{2k_{1}+2k_{2}})$ but $M$ is not diffeomorphic to $\C P^{\ell}\times S^{2k_{1}+2k_{2}}$ if $\rho \neq 0$ by Corollary~\ref{top classification of CASE2(1)}.
Therefore, we have the following lemma.
%%%%%%%%%%%%%%%%%%%%%%%%%%%%%%%%%%%%%%%%%%%%%%%%%%%%%%%%%%%%%%%%%%%%%%%%%%%
%Lemma 4
%%%%%%%%%%%%%%%%%%%%%%%%%%%%%%%%%%%%%%%%%%%%%%%%%%%%%%%%%%%%%%%%%%%%%%%%%%%
\begin{lemma}
\label{l>=4 trivial case}
    Let $M \in \mD$. Then $M$ is not cohomologically rigid but rigid by the cohomology ring and the Pontrjagin class in $\mM$, i.e., $\mD\subset \mR_{2}$.
\end{lemma}

%%%%%%%%%%%%%%%%%%%%%%%%%%%%%%%%%%%%%%%%%%%%%%%%%%%%%%%%%%%%%%%%%%%%%%%%%%%
%section 4.3.2
%%%%%%%%%%%%%%%%%%%%%%%%%%%%%%%%%%%%%%%%%%%%%%%%%%%%%%%%%%%%%%%%%%%%%%%%%%%
\subsubsection{The case when $M$ satisfies $H^{*}(M)\not\simeq H^{*}(\C P^{\ell}\times S^{2\deg w})$}
\label{sec:4.3.2}
Assume $H^{*}(M)\not\simeq H^{*}(\C P^{\ell}\times S^{2k_{1}+2k_{2}})$ and
the symbol $\mE$ represents the set of such manifolds.
By Section~\ref{sec:4.3.1}, we can easily show that $\mE$ is the set of manifolds
$B(\ell,\rho,k_{1},0)$ such that
\begin{enumerate}
\item $\rho\not=0$,
\item $k_{1}\le \ell$,
\item if $2k_{1}> \ell$ then $\rho\not\equiv_{2}0$.
\end{enumerate}
Note that the last statement is the same statement with that $\rho\not\equiv_{2}0$ or $2k_{1}\le\ell$.

Let $M_{1}$ and $M_{2}$ be elements in $\mE$ satisfying $H^\ast(M_1) \simeq H^\ast(M_2)$. Then we may assume that $M_{1}=B(\ell,\rho_{1},k_{1},0), M_{2}=B(\ell,\rho_{2},k_{1},0)$ and
\begin{align*}
H^{*}(M_{1})&=\Z[x_{1},w_{1}]/\langle x_{1}^{\ell+1}, w_{1}(w_{1}+\rho_{1}^{k_{1}}x_{1}^{k_{1}})\rangle \\
&\stackrel{\phi}{\simeq} \Z[x_{2},w_{2}]/\langle x_{2}^{\ell+1}, w_{2}(w_{2}+\rho_{2}^{k_{1}}x_{2}^{k_{1}})\rangle =H^{*}(M_{2}).
\end{align*}
With the method similar to that demonstrated in Section \ref{sec:4.3.1},
one can see that
if
$H^\ast(M_1) \stackrel{\phi}{\simeq} H^\ast(M_2)$
then
we have
\begin{eqnarray}
\label{phi-4.3.2}
\phi(x_{1})=\pm x_{2} \quad \text{ and } \quad \phi(w_{1})=a x_{2}^{k_{1}}+\epsilon w_{2},
\end{eqnarray}
for some $a \in \Z$ and $\epsilon=\pm 1$, and $\phi(w_{1}(w_{1}+(\rho_{1}x_{1})^{k_{1}})) = 0 \in H^\ast(M_2)$.
Now
\begin{eqnarray*}
& & \phi(w_{1}(w_{1}+(\rho_{1}x_{1})^{k_{1}})) \\
&=& (a x_{2}^{k_{1}}+\epsilon w_{2})^{2}+ (a x_{2}^{k_{1}}+\epsilon w_{2})(\pm\rho_{1}x_{2})^{k_{1}} \\
&=& a^{2}x_{2}^{2k_{1}}+2a\epsilon x_{2}^{k_{1}}w_{2}+w_{2}^{2}+a(\pm \rho_{1})^{k_{1}}x_{2}^{2k_{1}}+\epsilon(\pm \rho_{1})^{k_{1}}w_{2}x_{2}^{k_{1}} \\
&\equiv & a^{2}x_{2}^{2k_{1}}+2a\epsilon x_{2}^{k_{1}}w_{2}-w_{2}(\rho_{2}x_{2})^{k_{1}}+a(\pm \rho_{1})^{k_{1}}x_{2}^{2k_{1}}+\epsilon(\pm \rho_{1})^{k_{1}}w_{2}x_{2}^{k_{1}} \\
&=& ax_{2}^{2k_{1}}(a+(\pm \rho_{1})^{k_{1}})+\epsilon x_{2}^{k_{1}}w_{2}(2a-\epsilon\rho_{2}^{k_{1}}+(\pm \rho_{1})^{k_{1}}).
\end{eqnarray*}
Therefore,
we have that
if
$H^\ast(M_1) \simeq H^\ast(M_2)$
then
there is an $a \in \Z$ such that
\begin{eqnarray}
\label{eq4.3.2}
ax_{2}^{2k_{1}}(a+(\pm \rho_{1})^{k_{1}}) = x_{2}^{k_{1}}w_{2}(2a-\epsilon\rho_{2}^{k_{1}}+(\pm \rho_{1})^{k_{1}}) = 0 \in H^\ast(M_2).
\end{eqnarray}
On the other hand, if there is $a\in \Z$ which satisfies E.q.\ \eqref{eq4.3.2} then the above $\phi$ in \eqref{phi-4.3.2} provides the graded ring isomorphism $H^{*}(M_{1})\simeq H^{*}(M_{2})$.

If $2k_{1}\le \ell$, equivalently, $x_{2}^{2k_{1}} \neq 0$, then one can see that $|\rho_{2}|=|\rho_{1}|$ by \eqref{eq4.3.2}.
This implies that the vector bundles $(\gamma^{\otimes (-\rho_1)})^{\oplus k_{1}} \oplus \R$
and $(\gamma^{\otimes (-\rho_2)})^{\oplus k_{1}} \oplus \R$ are isomorphic as real vector bundles.
It follows that unit sphere bundles of these vector bundles are diffeomorphic.
Consequently, the manifold with $2k_{1}\le \ell$ is cohomologically rigid, i.e., the element in $\mR_{1}$.

Otherwise, i.e., $2k_{1}>\ell$, equivalently, $x_{2}^{2k_{1}} = 0$, then
$\rho_{1},\rho_{2}$ must be odd, i.e., $\rho_{1},\ \rho_{2}\not\equiv_{2}0$, by the definition of $\mE$.
Hence, by using \eqref{eq4.3.2},
if we put $a=\frac{\epsilon\rho_{2}^{k_{1}}+(\pm \rho_{1})^{k_{1}}}{2}$
then we get an isomorphism $H^{*}(M_{1})\stackrel{\phi}{\simeq}H^{*}(M_{2})$ for an
arbitrary odd numbers $\rho_1$ and $\rho_2$.
It follows from Corollary~\ref{top classification of CASE2(1)} that $M_1$ and $M_2$ are not necessarily diffeomorphic.
%Besides, for any $M_1$, there is $M_2$ which is not diffeomorphic to $M_1$ but has the same cohomology ring.
Thus, the manifold with $2k_{1}< \ell$ is not cohomologically rigid but rigid by the cohomology ring and the Pontrjagin class, i.e., the element in $\mR_{2}$.

In summary, the set $\mE$ can be divided into the following two sets:
\begin{eqnarray*}
\mE_{1}&=&\{B(\ell,\rho,k_{1},0) \ |\ \ell \ge4, \ 2<2k_{1}\le \ell,\ \rho\not=0  \}; \\
\mE_{2}&=&\{B(\ell,\rho,k_{1},0) \ |\ \ell \ge4, \ 2\le k_{1}\le \ell < 2k_{1},\ \rho\not\equiv_{2}0   \},
\end{eqnarray*}
and the following lemma holds.
%%%%%%%%%%%%%%%%%%%%%%%%%%%%%%%%%%%%%%%%%%%%%%%%%%%%%%%%%%%%%%%%%%%%%%%%%%%
%Lemma 5
%%%%%%%%%%%%%%%%%%%%%%%%%%%%%%%%%%%%%%%%%%%%%%%%%%%%%%%%%%%%%%%%%%%%%%%%%%%
\begin{lemma}
\label{l>=4 non-trivial case}
Let $M\in \mE\subset \mM$.
Then
\begin{enumerate}
\item $M$ is cohomologically rigid in $\mM$ if and only if $M\in \mE_{1}$, i.e., $\mE_{1}\subset \mR_{1}$;
\item $M$ is not cohomologically rigid but rigid by the cohomology ring and the Pontrjagin class in $\mM$ if and only if $M\in \mE_{2}$, i.e., $\mE_{2}\subset \mR_{2}$.
\end{enumerate}
\end{lemma}

The following proposition immediately follows from Lemmas \ref{l>=4 trivial case} and \ref{l>=4 non-trivial case}.
%%%%%%%%%%%%%%%%%%%%%%%%%%%%%%%%%%%%%%%%%%%%%%%%%%%%%%%%%%%%%%%%%%%%%%%%%%%
%Proposition CASE2 (1)
%%%%%%%%%%%%%%%%%%%%%%%%%%%%%%%%%%%%%%%%%%%%%%%%%%%%%%%%%%%%%%%%%%%%%%%%%%%
\begin{proposition}
\label{CASE2 (1)}
If $M\in \mM$ satisfies the condition of the CASE 2-(1), i.e., $\deg w>2$ (if and only if $2k_{1}+2k_{2} >2$) and $\ell \ge 4$, then
we can put $M=B(\ell,\rho,k_{1},k_{2})$ for $k_{1}>0$, $k_{2}\ge 0$ and
there are the following two cases:
\begin{enumerate}
\item $M\in \mR_{1}\subset \mM$ $\Leftrightarrow$ $M\in \mE_{1}$, i.e.,  $\rho\not=0$, $2k_{1}\le \ell$ and $k_{2}=0$;
\item $M\in \mR_{2}\subset \mM$ $\Leftrightarrow$ $M\in \mD\sqcup \mE_{2}$, i.e., $M$ satisfies one of the followings:
\begin{enumerate}
\item $M=B(\ell,0,k_{1},k_{2})\cong \C P^{\ell}\times S^{2k_{1}+2k_{2}}$;
\item $M=B(\ell,\rho,k_{1},k_{2})$ with $\rho\not=0$ and $k_{2}>0$;
%\item $M=B(\ell,\rho,k_{1},0)$ with $\rho\not=0$ and $k_{1}\ge \ell +1$;
%\item $M=B(\ell,\rho,k_{1},0)$ with $\rho\not=0$, $k_{1}<\ell+1\le 2k_{1}$.
\item $M=B(\ell,\rho,k_{1},0)$ with $\rho\not=0$, $\ell+1\le 2k_{1}$.
\end{enumerate}
\end{enumerate}
\end{proposition}

%%%%%%%%%%%%%%%%%%%%%%%%%%%%%%%%%%%%%%%%%%%%%%%%%%%%%%%%%%%%%%%%%%%%%%%%%%%
%section 4.4
%%%%%%%%%%%%%%%%%%%%%%%%%%%%%%%%%%%%%%%%%%%%%%%%%%%%%%%%%%%%%%%%%%%%%%%%%%%
\subsection{CASE 2 (2) : $\deg w >2$ and $\ell =2,\ 3$}
\label{sec:4.4}
Assume $\deg w >2$ and $\ell =2,\ 3$.
We first prove that this case is always rigid by the cohomology ring and the Pontrjagin class.

With a method similar to the one demonstrated in CASE 2 (1), we may put
\begin{eqnarray*}
M_1 = S((\gamma^{\otimes (-\rho_1)})^{\oplus k_{11}} \oplus \R^{2k_{12}+1}) = B(\ell,\rho_{1},k_{11},k_{12}) % S^{2\ell+1}\times_{S^{1}}S(\C_{\rho_1}^{ k_{11}} \oplus \R^{2k_{12}+1})
\end{eqnarray*}
and
\begin{eqnarray*}
M_2 = S((\gamma^{\otimes (-\rho_2)})^{\oplus k_{21}} \oplus \R^{2k_{22}+1}) = B(\ell,\rho_{2},k_{21},k_{22}) % S^{2\ell+1}\times_{S^{1}}S(\C_{\rho_2}^{ k_{21}} \oplus \R^{2k_{22}+1})
\end{eqnarray*}
for some $k_{11},\ k_{21}\in \N$ and $k_{12},\ k_{22}\ge 0$.
Let $\phi : H^\ast(M_1) \rightarrow H^\ast(M_2)$ be an isomorphism.
By \eqref{eq for x}, we have $\phi(x_{1})=\pm x_{2}$.
Note that, since $\ell=2$ or $3$, and from \eqref{p(M_i)}, we have that
$$
\phi( p(M_1)) = p(M_2) \iff k_{11} \rho_1^2 =k_{21} \rho_2^2
$$ because $p_j(M_i)$ vanishes for $j\ge 2$.

We remark that the second equation which holds in CASE 2 (1) does not hold in this case. Hence, we need to use another argument to show the rigidity by the cohomology and the Pontrjagin class; we will use a KO theoretical argument to do it.

Because of \cite{O}, in this case ($\ell=2,\ 3$), we have $KO(\C P^{\ell})\simeq \Z[y_{\ell}]/\langle y_{\ell}^{2}\rangle$,
where $y_{\ell}=r(\gamma)-2$ for the canonical line bundle $\gamma$ and the realification map
$r:K(\C P^{\ell})\to KO(\C P^{\ell})$.
Moreover, we have $r(\gamma^{\otimes n})=n^{2}y_{\ell}+2$ by \cite{O}.
Hence, for $i=1,~2$, we have that
\begin{eqnarray*}
r(\gamma^{\otimes (-\rho_{i})})=\rho_{i}^{2}r(\gamma)-2\rho_{i}^{2}+2.
\end{eqnarray*}
Hence, we have the following equation:
\begin{eqnarray*}
& &k_{11}r(\gamma^{\otimes (-\rho_{1})})+2k_{12}+1 \\
%=k_{11}(\rho_{1}^{2}y_{\ell}+2)+2k_{12}+1
&=&k_{11}\rho_{1}^{2}(r(\gamma)-2)+2(k_{11}+k_{12})+1 \\
&=&k_{21}\rho_{2}^{2}(r(\gamma)-2)+2(k_{21}+k_{22})+1 \\
%=k_{12}(\rho_{2}^{2}y_{\ell}+2)+2k_{22}+1
&=&k_{21}r(\gamma^{\otimes (-\rho_{2})})+2k_{22}+1
\end{eqnarray*}
in $KO(\C P^{\ell})$ by using $k_{11} \rho_1^2 =k_{21} \rho_2^2$ and $k_{11}+k_{12}=k_{21}+k_{22}$.
Therefore, we have that
\begin{eqnarray}
\label{condition of P II}
& & (\gamma^{\otimes (-\rho_{1})})^{\oplus k_{11}}\oplus \R^{2k_{12}+1}\equiv_{s} (\gamma^{\otimes (-\rho_{2})})^{\oplus k_{21}}\oplus \R^{2k_{22}+1} \nonumber \\
&\iff & k_{11} \rho_{1}^{2}=k_{21} \rho_{2}^{2}\quad {\rm and} \quad k_{11}+k_{12}=k_{21}+k_{22},
\end{eqnarray}
where $\eta\equiv_{s}\xi$ means two real vector bundles $\eta$ and $\xi$ that are stably isomorphic, i.e.,
there is a trivial vector bundle $\epsilon$ such that $\eta\oplus\epsilon\equiv \xi\oplus\epsilon$.

If $2k_{i1}+2k_{i2}+1>2\ell$ ($i=1,2$, and $\ell=2,3$), then these bundles are in the stable range, i.e., the dimension of fibre $2k_{i1}+2k_{i2}+1$ is strictly greater than that of the base space $\C P^{\ell}$; therefore, it follows from the stable range theorem (see e.g. \cite[Chapter 9]{H}) that the two bundles are isomorphic
$(\gamma^{\otimes (-\rho_{1})})^{\oplus k_{11}}\oplus \R^{2k_{12}+1}\equiv (\gamma^{\otimes (-\rho_{2})})^{\oplus k_{21}}\oplus \R^{2k_{22}+1}$.

Otherwise, i.e., $\ell=3$, $k_{i1}=k_{i2}=1$ ($i=1,2$), then we can easily show that $|\rho_{1}|=|\rho_{2}|$
since $k_{11}\rho_{1}^{2}=k_{21}\rho_{2}^{2}$; therefore, this case also satisfies that
$\gamma^{\otimes (-\rho_{1})}\oplus \R^{3}\equiv \gamma^{\otimes (-\rho_{2})}\oplus \R^{3}$.

It follows from the arguments above that if there is an isomorphism $\phi:H^{*}(M_{1})\to H^{*}(M_{2})$ such that $\phi$ preserves the Pontrjagin classes then
\begin{eqnarray*}
(\gamma^{\otimes (-\rho_{1})})^{\oplus k_{11}}\oplus \R^{2k_{12}+1}\equiv (\gamma^{\otimes (-\rho_{2})})^{\oplus k_{21}}\oplus \R^{2k_{22}+1};
\end{eqnarray*}
therefore, $M_{1}\cong M_{2}$.
Thus, we have the following lemma.
%%%%%%%%%%%%%%%%%%%%%%%%%%%%%%%%%%%%%%%%%%%%%%%%%%%%%%%%%%%%%%%%%%%%%%%%%%%
%Lemma 6
%%%%%%%%%%%%%%%%%%%%%%%%%%%%%%%%%%%%%%%%%%%%%%%%%%%%%%%%%%%%%%%%%%%%%%%%%%%
\begin{lemma}
\label{rigid by C-P II}
Assume $\deg w>2$ and $\ell= 2,\ 3$.
Then,
there is a graded ring isomorphism $\phi:H^{*}(M_{1})\to H^{*}(M_{2})$ such that $\phi(p(M_{1}))=p(M_{2})$
if and only if
$M_1$ and $M_2$ are diffeomorphic, i.e.,
this case is rigid by the cohomology ring and the Pontrjagin class.
\end{lemma}

We also have the following explicit topological classification of CASE 2 (2) by using \eqref{condition of P II} and Lemma~\ref{rigid by C-P II}.
%%%%%%%%%%%%%%%%%%%%%%%%%%%%%%%%%%%%%%%%%%%%%%%%%%%%%%%%%%%%%%%%%%%%%%%%%%%
%Corollary 4
%%%%%%%%%%%%%%%%%%%%%%%%%%%%%%%%%%%%%%%%%%%%%%%%%%%%%%%%%%%%%%%%%%%%%%%%%%%
\begin{corollary}
\label{top classification of CASE2(2)}
Assume $\deg w>2$ and $\ell=2,\ 3$.
Then, $B(\ell,\rho,k_{1},k_{2})\cong B(\ell',\rho',k_{1}',k_{2}')$ if and only if
$\ell=\ell'$, $k_{1}+k_{2}=k_{1}'+k_{2}'$ and
$k_{1} \rho^{2}=k_{1}' (\rho')^{2}$.
\end{corollary}

Here we exhibit some non-trivial examples.
\begin{example} \label{example:counterexample}
The following two manifolds are diffeomorphic because $H^{*}(M_{1})\simeq H^{*}(M_{2})$ and $p_{1}(M_{1})=4x_{1}^{2}$ and
$p_{1}(M_{2})=4x_{2}^{2}$ ($x_{i}\in H^{2}(M_{i})$):
\begin{eqnarray*}
M_{1}=S^{7}\times_{S^{1}}S(\C_{2}^{1}\oplus \R^{9});\quad
M_{2}=S^{7}\times_{S^{1}}S(\C_{1}^{4}\oplus \R^{1}).
\end{eqnarray*}
The following manifold has the same cohomology ring as the above two manifolds,
but this manifold is not diffeomorphic to the above manifolds
because $p_{1}(M)=16x^{2}$ for $x\in H^{2}(M)$:
\begin{eqnarray*}
M=S^{7}\times_{S^{1}}S(\C_{2}^{4}\oplus \R^{1}). \\
\end{eqnarray*}
\end{example}

With the method similar to that demonstrated in Sections \ref{sec:4.3.1} and \ref{sec:4.3.2},
we have
$M\in\mR_{1}\subset \mM$ if and only if $\rho\not=0$, $2k_{1}\le \ell$ and $k_{2}=0$.
However, it follows from $\ell=2$ or $3$ that $k_{1}=1$.
This gives a contradiction to $(k_{1},k_{2})\not=(1,0)$ (see Section \ref{sec:4.2}).
Therefore, we have the following proposition by the method similar to that demonstrated in Sections \ref{sec:4.3.1} and \ref{sec:4.3.2}.
%%%%%%%%%%%%%%%%%%%%%%%%%%%%%%%%%%%%%%%%%%%%%%%%%%%%%%%%%%%%%%%%%%%%%%%%%%%
%Proposition CASE2 (2)
%%%%%%%%%%%%%%%%%%%%%%%%%%%%%%%%%%%%%%%%%%%%%%%%%%%%%%%%%%%%%%%%%%%%%%%%%%%
\begin{proposition}
\label{CASE2-(2)}
If $M\in \mM$ satisfies the condition of the CASE 2-(2), i.e., $\deg w>2$ (if and only if $2k_{1}+2k_{2} >2$) and $\ell=2,\ 3$, then
we can put $M=B(\ell,\rho,k_{1},k_{2})$ for $k_{1}>0$, $k_{2}\ge 0$ and $M\in\mR_{2}\subset \mM$.
Furthermore, $M$ satisfies one of the followings:
\begin{enumerate}
\item $M=B(\ell,0,k_{1},k_{2})\cong \C P^{\ell}\times S^{2k_{1}+2k_{2}}$;
\item $M=B(\ell,\rho,k_{1},k_{2})$ with $\rho\not=0$ and $k_{2}>0$;
\item $M=B(\ell,\rho,k_{1},0)$ with $\rho\not=0$, $\ell+1\le 2k_{1}$ (i.e., $k_{1}\ge 2$).
\end{enumerate}
\end{proposition}

%%%%%%%%%%%%%%%%%%%%%%%%%%%%%%%%%%%%%%%%%%%%%%%%%%%%%%%%%%%%%%%%%%%%%%%%%%%
%section 4.5
%%%%%%%%%%%%%%%%%%%%%%%%%%%%%%%%%%%%%%%%%%%%%%%%%%%%%%%%%%%%%%%%%%%%%%%%%%%
\subsection{CASE 2 (3) : $\deg w >2$ and $\ell=1$}
\label{sec:4.5}

Assume $\deg w >2$ and $\ell=1$.
In this case,
one can easily see that $H^{*}(M)\simeq H^{*}(S^{2k_{1}+2k_{2}}\times S^{2})$ from \eqref{cohomology} and \eqref{f_i}.
Moreover,
we may put
\begin{eqnarray*}
M=S^{3}\times_{S^{1}}S(\C_{\rho}^{k_{1}}\oplus\R^{2k_{2}+1})=B(1,\rho,k_{1},k_{2})
\end{eqnarray*}
for $\rho\in \Z$, $k_{1}>0$ and $k_{2}\ge 0$, where $(k_{1},k_{2})\not =(1,0)$ and $\deg w=2k_{1}+2k_{2}$.
Therefore, if $H^{*}(B(1,\rho,k_{1},k_{2}))\simeq H^{*}(B(1,\rho',k_{1}',k_{2}'))$ then we have $k_{1}+k_{2}=k_{1}'+k_{2}'$ only.

Using Propositions~\ref{prop2.2} and~\ref{prop2.3}, we have
\begin{eqnarray}
\label{SW-class in III}
p(M)=1 \quad \text{ and } \quad w(M)=1+k_{1}\rho x.
\end{eqnarray}
It follows that the Pontrjagin class does not distinguish diffeomorphism types of this case.

Recall that $S^{n-1}$-bundles over $S^{2}$ are classified by
continuous maps from $S^{2}$ to $G_{n}=B O(n)$ up to homotopy and
$\pi_2(G_n) \simeq \Z_2$ for $n>2$ (see e.g. \cite{S}).
We can easily show that this $\Z_{2}$ is generated by $w_{2}(M)$, i.e., $k_{1}\rho x$ in our case (see \eqref{SW-class in III}).
If we fix the dimension of its fibre,
it follows that there are just two sphere bundles over $S^{2}$, i.e., the trivial bundle and the non-trivial bundle.
In our case, if $k_{1}$ or $\rho$ is even, then $S^{3}\times_{S^{1}}(\C_{\rho}^{k_{1}}\oplus \R^{2k_{2}+1})$ is always trivial bundle because its Stiefel-Whitney class is trivial.
Hence, if $k_{1}$ or $\rho$ is even, then $M\cong S^{2k_{1}+2k_{2}}\times S^{2}$.
If $k_{1}$ and $\rho$ are odd, then
$S^{3}\times_{S^{1}}(\C_{\rho}^{k_{1}}\oplus \R^{2k_{2}+1})$ is the non-trivial bundle.

Therefore, the following proposition holds:
\begin{proposition}
\label{CASE2-(3)}
If the CASE 2-(3) holds, i.e., $\ell=1$ and $2k_{1}+2k_{2} >2$, then
$M=B(1,\rho,k_{1},k_{2})$ for $k_{1}>0$, $k_{2}\ge 0$ and
$M$ is not rigid by the cohomology ring and the Pontrjagin class but rigid by the cohomology ring and the Stiefel-Whitney class, i.e.,
$M\in \mR_{3}$.
\end{proposition}

We also have the following explicit topological classification of CASE 2 (3).
%%%%%%%%%%%%%%%%%%%%%%%%%%%%%%%%%%%%%%%%%%%%%%%%%%%%%%%%%%%%%%%%%%%%%%%%%%%
%Corollary 4
%%%%%%%%%%%%%%%%%%%%%%%%%%%%%%%%%%%%%%%%%%%%%%%%%%%%%%%%%%%%%%%%%%%%%%%%%%%
\begin{corollary}
\label{top classification of CASE2(3)}
If the CASE 2-(3) holds, i.e., $\ell=1$ and $2k_{1}+2k_{2} >2$, then
their diffeomorphism types are classified as follows:
\begin{enumerate}
\item $B(1,\rho,k_{1},k_{2})\cong B(1,0,k_{1},k_{2})=S^{2}\times S^{2k_{1}+2k_{2}}$ if and only if $k_{1}\equiv_{2}0$ or
$\rho\equiv_{2}0$;
\item $B(1,\rho,k_{1},k_{2})\cong B(1,1,k_{1},k_{2})=S^{3}\times_{S^{1}} S(\C^{k_{1}}_{1}\oplus\R^{2k_{2}+1})$ if and only if otherwise, i.e., $k_{1}\equiv_{2}1$ and $\rho\equiv_{2}1$.
\end{enumerate}
\end{corollary}

By Propositions~\ref{CASE1},~\ref{CASE2 (1)},~\ref{CASE2-(2)} and~\ref{CASE2-(3)}, we have Theorem \ref{thm-main}.

%%%%%%%%%%%%%%%%%%%%%%%%%%%%%%%%%%%%%%%%%%%%%%%%%%%%%%%%%%%%%%%%%%%%%%%%%%%
%section 5
%%%%%%%%%%%%%%%%%%%%%%%%%%%%%%%%%%%%%%%%%%%%%%%%%%%%%%%%%%%%%%%%%%%%%%%%%%%

\section{Further Study}
\label{sec:5}
A homotopy cell is called a \emph{homotopy polytope} if any multiple intersection of faces is connected whenever it is non-empty. We note that the set of homotopy polytopes contains the set of simple polytopes while it is included in the set of homotopy cells. A torus manifold with
a
homotopy polytope
as its
orbit space is also an interesting object in toric topology. Indeed, Masuda and Suh \cite{MS} expected that toric theory can be applied to the families of such manifolds in a nice way, so they asked the cohomological rigidity problem for the above two classes.

By Theorem~\ref{thm-main}, we answered negatively to the cohomological rigidity problem for the family of torus manifolds whose orbit spaces are homotopy cells. However, we still do not know the answer for the case
where the orbit space is a homotopy polytope.

\begin{problem}
Are torus manifolds $X$ and $Y$ whose orbit spaces are homotopy polytopes diffeomorphic (or homeomorphic) if $H^\ast(X) \cong H^\ast(Y)$ as graded rings?
\end{problem}

Moreover, we may ask the following question from our result.

\begin{problem}
Are two torus manifold with homotopy cell as its orbit spaces (or codimension one extended actions) are classified by their cohomology ring and real characteristic classes?
\end{problem}

\section*{Acknowledgements}
Finally the authors would like to thank Professors Zhi L\"u and Dong Youp Suh for providing excellent circumstances to do research.
They also would like to thank for
Professors Mikiya Masuda and Andreas Holmsen for their invaluable advices and comments.

\end{document}